\documentclass[10pt]{amsart}

\newcommand{\comm}[1]{}

\numberwithin{equation}{section}
\let\newpf\proof \let\proof\relax 
\newenvironment{pf}{\newpf[\proofname]}{\qed\endtrivlist}

\newcommand{\LE}{\mathcal{E}}
\renewcommand{\P}{\mathbb{P}}
\newcommand{\RR}{\mathfrak{R}}
\renewcommand{\l}{{\underline{l}}}

\def\be{\begin{equation}}
\def\ee{\end{equation}}

\def\ba{{\begin{align}}}
\def\ea{{\end{align}}}

\newcommand\ssigma{\mathfrak{S}^0}
\newcommand\sssigma{\mathfrak{S}}

\newcommand{\rank}{\operatorname{rank}}
\newcommand{\Ext}{{\mathrm {\Lambda}}}
\newcommand{\grass}{{\mathrm {Grass}}}
\newcommand{\iso}{{\mathrm {Iso}}}
\newcommand{\GL}{{\mathrm {GL}}}
\newcommand{\SL}{{\mathrm {SL}}}

\newtheorem{main}{Theorem}
\newtheorem{thm}{Theorem}[section]
\newtheorem{cor}[thm]{Corollary}

\newtheorem{lemma}[thm]{Lemma}

\theoremstyle{remark}
\newtheorem{rem}[thm]{Remark}

\newtheorem{problem}[thm]{Problem}
\theoremstyle{definition}
\newtheorem{definition}[thm]{Definition}

\newcommand{\dist}{\operatorname{dist}}
\newcommand{\id}{\operatorname{id}}
\newcommand{\Ker}{\operatorname{Ker}}

\renewcommand{\AA}{{\mathcal A}}
\newcommand{\BB}{{\mathcal B}}
\newcommand{\CC}{{\mathcal C}}
\newcommand{\DD}{{\mathcal D}}
\newcommand{\EE}{{\mathcal E}}

\newcommand{\FF}{{\mathcal F}}

\newcommand{\LL}{{\mathcal L}}
\newcommand{\MM}{{\mathcal M}}
\newcommand{\NN}{{\mathcal N}}

\newcommand{\XX}{{\mathcal X}}
\newcommand{\YY}{{\mathcal Y}}

\newcommand{\N}{{\mathbb N}}

\newcommand{\R}{{\mathbb R}}

\newcommand{\Z}{{\mathbb Z}}

\begin{document}

\title[Proof of the Zorich-Kontsevich conjecture]
{Simplicity of Lyapunov spectra:\\
Proof of the Zorich-Kontsevich Conjecture}

\author{Artur Avila and Marcelo Viana}

\thanks{Work supported by the Brazil-France Agreement in
Mathematics. M.V. is also supported by Pronex and Faperj.}

\address{
CNRS UMR 7599, Laboratoire de Probabilit\'es et Mod\`eles
Al\'eatoires, Universit\'e Pierre et Marie Curie, Bo\^\i te
Postale 188, 75252 Paris Cedex 05, France }
\email{artur@ccr.jussieu.fr}

\address{
IMPA -- Estr. D. Castorina 110, Jardim Bot\^anico,
22460-320 Rio de Janeiro -- Brazil.
}

\email{viana@impa.br}

\begin{abstract}
We prove the Zorich-Kontsevich conjecture that the non-trivial
Lyapunov exponents of the Teichm\"uller flow on (any connected
component of a stratum of) the moduli space of Abelian
differentials on compact Riemann surfaces are all distinct. By
previous work of Zorich and Kontsevich, this implies the existence
of the complete asymptotic Lagrangian flag describing the behavior
in homology of the vertical foliation in a typical translation
surface.

\end{abstract}

\maketitle

\section{Introduction}

A {\it translation surface} is a compact orientable surface
endowed with a flat metric with finitely many conical
singularities and a unit parallel vector field. Another equivalent
description is through complex analysis: a translation surface is
a pair $(M,\omega)$ where $M$ is a Riemann surface and $\omega$ is
an  Abelian differential, that is, a complex holomorphic $1$-form:
the metric is then given by $|\omega|$ and a unit parallel vector
field is specified by $\omega\cdot v=i$. We call {\it vertical
flow} the flow $\phi_t$ along the vector field $v$. It is defined
for all times except for finitely many orbits that meet the
singularities in finite time.

The space of all translation surfaces of genus $g$, modulo isometries
preserving the parallel vector fields, is thus identified with the
moduli space of Abelian differentials $\MM_g$.
If we also specify the orders $\kappa$ of the zeroes of $\omega$, the space
of translation surfaces, modulo isometry, is identified with a {\it stratum}
$\MM_{g,\kappa} \subset \MM_g$.
Each such stratum is a complex affine variety, and so is endowed with a
Lebesgue measure class.

\subsection{The Zorich phenomenon}

Pick a typical (with respect to the Lebesgue measure class) translation surface
$(M,\omega)$ in $\MM_{g,\kappa}$. Consider segments of orbit
$\Gamma^T_x=\{\phi_t(x),\, 0 \leq t \leq e^T\}$, where $x$ is chosen
arbitrarily such that the vertical flow is defined for all times,
and let $\gamma^T_x$ be a closed loop obtained by concatenating $\Gamma^T_x$
with a segment  of smallest possible length joining $\phi_{e^T}(x)$ to $x$.
Then $[\gamma^T_x] \in H_1(M)$ is asymptotic to $e^T c$, where
$c \in H_1(M) \setminus \{0\}$ is the Schwartzman  \cite {S} {\it asymptotic
cycle}.

If $g=1$ then this approximation is quite good: the deviation from the line
$F_1$ spanned by the asymptotic cycle is bounded. When $g=2$, one gets a richer
picture: $[\gamma^T_x]$ oscillates around $F_1$ with amplitude roughly
$e^{\lambda_2 T}$, where $0<\lambda_2<1$.  Moreover, there is an asymptotic
isotropic $2$-plane $F_2$: deviations from $F_2$ are bounded.
More generally, in genus $g>1$, there exists an asymptotic Lagrangian flag,
that is, a sequence of nested isotropic spaces $F_i$, $1 \leq i \leq g$
of dimension $i$ and numbers $1>\lambda_2>\cdots>\lambda_g>0$ such that
$[\gamma^T_x]$ oscillates around $F_i$ with amplitude roughly
$e^{\lambda_{i+1}T}$, $1 \leq i \leq g-1$ and the deviation from $F_g$ is
bounded:
\be
\limsup_{T\to\infty}  \frac {1} {T}
\ln \dist([\gamma^T_x],F_i)=\lambda_{i+1}
 \quad\text{for every } 1 \leq i \leq g-1
\ee
\be
\sup \dist([\gamma^T_x],F_g)<\infty.
\ee
Moreover, the deviation spectrum $\lambda_2>\cdots>\lambda_g$ is universal:
it depends only on the connected component of the stratum to which $(M,\omega)$
belongs.

The picture we just described is called the Zorich phenomenon, and was
discovered empirically by Zorich \cite {Z1}.
It was shown by Kontsevich and Zorich  \cite {Z1,Z2,KZ} that this picture
would follow from a statement about the Lyapunov exponents of the
Kontsevich-Zorich cocycle, that we now discuss.

\subsection{The Kontsevich-Zorich cocycle}

The Kontsevich-Zorich cocycle can be described roughly as follows.
It is better to work on (the cotangent bundle of) the
Teichm\"uller space, that is, to consider translation structures
on a surface up to isometry isotopic to the identity. The moduli
space is obtained by taking the quotient by the modular group. Let
$(M_T,\omega_T)$ be obtained from $(M,\omega)$ by applying the
Teichm\"uller flow for time $T$, that is, we change the flat
metric by contracting (by $e^{-T}$) the vertical direction, and
expanding (by $e^T$) the orthogonal horizontal direction. The
Riemann surface $M_T$ now looks very distorted but it can be
brought back to a fundamental domain in the Teichm\"uller space by
applying an element $\Phi_T(M,\omega)$ of the mapping class group.
The action of $\Phi_T$ on the cohomology $H^1(M)$ is, essentially,
the Kontsevich-Zorich cocycle.

The Kontsevich-Zorich cocycle is, thus, a linear cocycle over the
Teichm\"uller flow. It was shown by Masur \cite {M} and Veech
\cite {V1} that the Teichm\"uller flow restricted to each
connected component of strata is ergodic, provided we  normalize
the area. The cocycle is measurable (integrable), and so it has
Lyapunov exponents $\lambda_1 \geq \cdots \geq \lambda_{2g}$. The
statement about the cocycle that implies the Zorich phenomenon is
that the cocycle has simple Lyapunov spectrum, that is, its
Lyapunov exponents are all distinct: \be\label{simples1}
\lambda_1>\lambda_2>\cdots>\lambda_{2g-1}>\lambda_{2g} \ee (see
\cite{Z4}, Theorem 2 and Conditional Theorem 4). The $\lambda_i$,
$1<i \leq g$ are the same that appear in the description of the
Zorich phenomenon, and the asymptotic flag is related to the
Oseledets decomposition. One has $\lambda_i=-\lambda_{2g-i+1}$ for
all $i$, because the cocycle is symplectic. It is easy to see that
$\lambda_1=1=-\lambda_{2g}$. The $\lambda_i$ are also related to
the Lyapunov exponents of the Teichm\"uller flow on the
corresponding connected component of strata: the latter are,
exactly, \be\label{simples2} 2  \geq  1+\lambda_2 \geq \cdots \geq
1+\lambda_g \geq 1 = \cdots = 1 \geq 1+\lambda_{g+1}  \geq \cdots
\geq 1+\lambda_{2g-1} \geq 0 \ee together with their symmetric
values. Zero is a simple exponent, corresponding to the flow
direction. There are $\sigma-1$ exponents equal  to $1$ (and to
$-1$) arising from the action on relative cycles joining the
$\sigma$ singularities. It is clear that \eqref{simples1} is
equivalent to saying that all the inequalities in \eqref{simples2}
are strict.

It follows from the work of Veech \cite {V2} that $\lambda_2<1$,
and so the first and the last inequalities in \eqref{simples2} are
strict. Thus, the Teichm\"uller flow is non-uniformly hyperbolic.
Together with Kontsevich \cite {K}, Zorich conjectured that the
$\lambda_i$ are, indeed, all distinct. They also established
formulas for the sums of Lyapunov exponents, but it has not been
possible to use them for proving this conjecture. The fundamental
work of Forni \cite {F} established that
$\lambda_g>\lambda_{g+1}$, which is the same as saying that no
exponent vanishes. This means that $0$ does not belong to the
Lyapunov spectrum, that is, the Kontsevich-Zorich cocycle is
non-uniformly hyperbolic. Besides giving substantial information
on the general case, this result settles the conjecture in the
particular case $g=2$, and has also been used to obtain other
dynamical properties of translation flows \cite {AF}.

\subsection{Main result}

Previously to the introduction of the Kontsevich-Zorich cocycle,
Zorich had already identified and studied a discrete time version,
called the Zorich cocycle. Its precise definition will be recalled
in section \ref {rauzy}. While the base dynamics of the
Kontsevich-Zorich cocycle is the Teichm\"uller flow on the space
of translation surfaces, the basis of the Zorich cocycle is a
renormalization dynamics in the space of interval exchange
transformations. The link between interval exchange
transformations and translation flows is well known: the first
return map to a horizontal cross-section to the vertical flow is
an interval exchange transformation, and the translation flow can
be reconstructed as a special suspension over this one-dimensional
transformation. Typical translation surfaces in the same connected
component of strata are special suspensions over interval exchange
transformations whose combinatorics belong to the same {\it Rauzy
class}. The Zorich cocycle is a measurable (integrable) linear
cocycle (over the renormalization dynamics) acting on a vector
space $H$ that can be naturally identified with $H_1(M)$. So, it
has Lyapunov exponents $\theta_1 \geq \cdots \geq \theta_{2g}$.
The $\theta_i$ are linked to the $\lambda_i$ by
$$
\lambda_i=\frac {\theta_i} {\theta_1} \quad\text{for all } 1 \leq i \leq 2g,
$$
so that properties of the $\theta_i$ can be deduced from those of
the $\lambda_i$, and conversely. The Zorich Conjecture (see \cite
{Z3}, Conjecture 1 or \cite {Z4}, Conjecture 2) states that the
Lyapunov spectrum of the Zorich cocycle is simple, that is, all
the $\theta_i$ are distinct. By the previous discussion, it
implies the full picture of the Zorich phenomenon. Here we prove
this conjecture:

\begin{main}\label{main}
The Zorich cocycle has simple Lyapunov spectrum on every Rauzy class.
\end{main}

The most important progress in this direction so far, the work of
Forni \cite {F}, was via the Kontsevich-Zorich cocycle. Here we
address the Zorich cocycle directly. Though many ideas can be
formulated in terms of the Kontsevich-Zorich cocycle, our approach
is mostly dynamical and does not involve the extra geometrical and
complex analytic structures present in the Teichm\"uller flow
(particularly the $\SL(2,\R)$ action on the moduli space) that
were crucial in \cite {F}. Our arguments contain, in particular, a
new proof of Forni's main result.

A somewhat more extended discussion of the Zorich-Kontsevich
conjecture, including an announcement of our present results, can
be found in \cite{AV0}.

\begin{rem}
Besides the Zorich phenomenon, the Lyapunov exponents of the
Zorich cocycle are also linked to the behavior of ergodic averages
of interval exchange transformations (a first result in this
direction is given in \cite {Z2}) and translation flows or
area-preserving flows on surfaces \cite {F}. Let us point out that
it is now possible to treat the case of interval exchange
transformations in a very elegant way, using the results of
Marmi-Moussa-Yoccoz \cite {MMY}. A result for translation flows
can then be recovered by suspension, which also implies the result
for area-preserving flows.
\end{rem}

\subsection{A geometric motivation for the proof}

Our proof of the Zorich conjecture has two distinct parts:
\begin{enumerate}
\item A general criterion for the simplicity of the Lyapunov spectrum of
locally constant cocycles.
\item A combinatorial analysis of Rauzy diagrams to show that the criterion
can be applied to the Zorich cocycle on any Rauzy class.
\end{enumerate}
The basic idea of the criterion is that it suffices to find orbits of the
base dynamics over which the cocycle exhibits certain forms of behavior,
that we call ``twisting'' and ``pinching''.
Roughly speaking, twisting means that certain families of subspaces are put
in general position, and pinching means that a large part of the Grassmannian
(consisting of subspaces in general position) is concentrated in a small
region, by the action of the cocycle. We will be a bit more precise in a while.
We call the cocycle {\it simple} if it meets both requirements (notice that
``simple'' really means that the cocycle's behavior is quite rich).
Then, according to our criterion, the Lyapunov spectrum is simple.

It should be noted that the orbits on which these types of
behavior are observed are very particular and, a priori,
correspond to zero measure subsets of the orbits. Nevertheless,
they are able to ``persuade'' almost every orbit to have a simple
Lyapunov spectrum. For this, one assumes that the base dynamics is
rather chaotic (which is the case for the Teichm\"uller flow). In
a purely random situation, this persuasion mechanism has been
understood for quite some time, through the works of Furstenberg
\cite{Fu}, Guivarc'h and Raugi \cite{GR}, and Gol'dsheid and
Margulis \cite {GM}.  That this happens also in chaotic, but not
random, situations was unveiled by the works of Ledrappier
\cite{L} and Bonatti, Gomez-Mont, Viana \cite {BGMV,BV,V}. The
particular situation needed for our arguments is, however, not
covered in those works, and was dealt with in our twin paper
\cite{AV}.

The present paper is devoted to part (2) of the proof. However, to
make the paper self-contained, in Appendix~\ref{simplicity
criterion} we also prove an instance of the criterion (1) that
covers the present situation.

The basic idea to prove that the Zorich cocycle is ``rich'', in
the sense described above, is to use induction on the complexity.
The geometric motivation is more transparent when one thinks in
terms of the Kontsevich-Zorich cocycle. As explained, we want to
find orbits of the Teichm\"uller flow inside any connected
component of a stratum $\CC$ with some given behavior. To this
end, we look at orbits that spend a long time near the boundary of
$\CC$. While there, these orbits pick up the behavior of the
boundary dynamics of the Teichm\"uller flow, which contains the
dynamics of the Teichm\"uller flow restricted to connected
components of strata $\CC'$ with simpler combinatorics
(corresponding to certain ways to degenerate $\CC$).

This is easy to make sense of when the stratum $\CC$ is not closed
in the moduli space, since the whole Teichm\"uller flow provides a
broader ambient dynamics where everything takes place. It is less
clear how to formalize the idea when $\CC$ is closed (this is the
case when there is only one singularity). In this case, the
boundary dynamics corresponds to the Teichm\"uller flow acting on
surfaces of smaller genuses, and here we will not attempt to
describe geometrically what happens. Let us point out that
Kontsevich-Zorich~\cite {KZ} considered the inverse of such a
degeneration process, that they called ``bubbling a handle''. It
is worth noting that, though they used many different techniques,
and as much geometric reasoning as possible, this was one of the
few steps that needed to be done using combinatorics of interval
exchange transformations (encoded in Rauzy diagrams).

On the other hand, our degeneration process is very simple when
viewed in terms of interval exchange transformations: we just make
one interval very small. This small interval remains untouched by
the renormalization process for a very long time, while the other
intervals are acted upon by a degenerate renormalization process.
This is what allows us to put in place an inductive argument. To
really control the effect, we must choose our small interval very
carefully. It is also sometimes useful to choose particular
permutations in the Rauzy class we are analyzing. A particularly
sophisticated choice is needed when we must change the genus of
the underlying translation surface: at this point we use Lemma 20
of Kontsevich-Zorich~\cite{KZ}, which allowed them to obtain the
inverse process of ``bubbling a handle''.

Let us say a few more words on how the Zorich cocycle will be
shown to be simple. The fact that the cocycle is symplectic is
important for the arguments. By induction, we show that it acts
minimally on the space of Lagrangian flags. This is used to derive
that the cocycle is ``twisting'': certain orbits can be used to
put families of subspaces in general position. In the induction,
there must be some gain of information at each step, when we must
change genus: in this case, this gain regards the action of the
cocycle on lines, and it comes from the rather easy fact that this
action is minimal.

Also by induction, we show that certain orbits of the cocycle are
``pinching'': they take a large amount of the Grassmannian and
concentrate it into a small region. Here the gain of information
when we must change genus has to come from the action on
Lagrangian spaces, and it is far from obvious. One can use Forni's
theorem \cite {F} and, indeed, we did so in a previous version of
the arguments. However, the proof we will give is independent of
his result, so that our work gives a new proof of Forni's main
theorem. Indeed, in the (combinatorial) argument that we will
give, the pinching of Lagrangian subspaces comes from orbits that
have a pair of zero Lyapunov exponents, but present some parabolic
behavior in the central subspace.

\subsection{Outline of the paper}

Section \ref {gb} gives general background on linear and
symplectic actions of monoids. Section \ref {rauzy} collects
well-known material on interval exchange transformations. We
introduce the Zorich cocycle and discuss the combinatorics of the
Rauzy diagram. The presentation follows \cite {MMY} closely. The
matrices appearing in the Zorich cocycle can be studied in terms
of the natural symplectic action of a combinatorial object that we
call the Rauzy monoid. We give properties of the projective action
of the Rauzy monoids, and discuss some special elements in Rauzy
classes. In section \ref {4}, we introduce the twisting and
pinching properties of symplectic actions of monoids, and use them
to define the notion of simplicity for monoid actions, which is
our basic sufficient condition for simplicity of the Lyapunov
spectrum. In sections \ref {5} and \ref {pinra} we prove the
twisting and pinching properties, and thus simplicity, for the
action of the Rauzy monoid. The proof involves a combinatorial
analysis of relations between different Rauzy classes, as
explained before. In section \ref {simpl}, we state the basic
sufficient criterion for simplicity of the Lyapunov spectrum
(Theorem \ref{simplicity criterion}) and show that it is satisfied
by the Zorich cocycle on any Rauzy class. Theorem  \ref {main}
follows. For completeness, in  Appendix~\ref{simplicity criterion}
we give a proof of that sufficient criterion.

\medskip

{\bf Acknowledgments:} We would like to thank Jean-Christophe
Yoccoz for several inspiring discussions, through which he
explained to us his view of the combinatorics of interval exchange
transformations.  We also thank him, Alexander Arbieto, Giovanni
Forni, Carlos Matheus, and Weixiao Shen for listening to many
sketchy ideas while this work developed, and Evgeny Verbitsky for
a discussion relevant to Section~\ref{simpl}.

\section{General background} \label {gb}

\subsection{Grassmannian structures}

Throughout this paper, all vector spaces are finite dimensional
vector spaces over $\R$. The notation $\P H$ always represents the
projective space, that is, the space of lines ($1$-dimensional
subspaces) of a vector space $H$. More generally,  $\grass(k,H)$
will represent the {\it Grassmannian} of $k$-planes, $1 \leq k
\leq \dim H-1$ in the space $H$. For $0 \leq k \leq \dim H$, we
denote by $\Ext^k(H)$ the $k$-th exterior product of $H$. If $F
\in \grass(k,H)$ is spanned by linearly independent vectors
$v_1$,\ldots,$v_k$ then the exterior product of the $v_i$ is
defined up to multiplication by a scalar. This defines an
embedding $\grass(k,H) \to\P\Ext^k(H)$.

A {\it geometric line} in $\Ext^k H$ is a line that  is contained
in $\grass(k,H)$. The duality between $\Ext^k H$ and $\Ext^{\dim
H-k} H$ allows one to define a {\it geometric hyperplane} of
$\Ext^k H$ as the dual to a geometric line of $\Ext^{\dim H-k} H$.
A {\it hyperplane section} is the intersection with $\grass(k,H)$
of the projectivization of a geometric hyperplane in $\Ext^k H$.
In other words, it is the set of all $F \in \grass(k,H)$ having
non-trivial intersection with a given $E \in\grass(\dim H-k,H)$.
Thus, hyperplane sections are closed subsets of $\grass(k,H)$ with
empty interior.  In particular, $\grass(k,H)$ can not be written
as a finite, or even countable, union of hyperplane sections.

We will call {\it linear arrangement}  in $\Ext^k H$ any finite union of finite
intersections of  geometric hyperplanes. The intersection of the
projectivization of a linear arrangement with $\grass(k,H)$ will be called a
linear arrangement in $\grass(k,H)$.  It will be called non-trivial if it is
neither empty nor the whole $\grass(k,H)$.

The {\it flag space} $\FF(H)$ is the set of all $(F_i)_{i=1}^{\dim H-1}$ where
$F_i$ is a subspace of $H$ of dimension $i$ and $F_i \subset F_{i+1}$ for all
$i$. It is useful to see $\FF(H)$ as a fiber bundle $\FF(H) \to \P H$,
through the projection $(F_i)_{i=1}^{\dim H-1} \mapsto F_1$.
The fiber over $\lambda \in \P H$ is naturally isomorphic to $\FF(H/\lambda)$,
via the isomorphism
\be
(F_i)_{i=1}^{\dim H-1} \mapsto (F_{i+1}/\lambda)_{i=1}^{\dim H-2}.
\ee

We are going to state a few simple facts about the family of
finite non-empty unions of linear subspaces of any given vector
space, and deduce corresponding statements for the family of
linear arrangements. The proofs are elementary, and we leave them
for the reader.

\begin{lemma}\label{l.linear1}

Any finite union $W$ of linear subspaces of some vector space
admits a canonical expression $W=\cup_{V\in\XX} V$ which is
minimal in the following sense: if $W=\cup_{V\in\YY} V$ is another
way to express $W$ as a union of linear subspaces then
$\XX\subset\YY$.

\end{lemma}

It is clear that the family of finite unions of linear subspaces is closed
under finite unions and finite intersections. The next statement implies that
it is even closed under arbitrary intersections.

\begin{lemma}\label{l.linear3}

If $\{V^\alpha: \alpha\in A\}$ is an arbitrary family of finite
unions of linear subspaces of some vector space, then
$\cap_{\alpha\in A} V^\alpha$ coincides with the intersection of
the $V^\alpha$ over a finite subfamily.

\end{lemma}

%

\begin{cor}\label{c.linear4}

Any totally ordered (under inclusion) family $\{V^\alpha:
\alpha\in A\}$ of finite unions of linear subspaces of some vector
space is well ordered.

\end{cor}

%

\begin{cor}\label{c.linear5}

If $V$ is a finite union of linear subspaces of some vector space
and $x$ is an isomorphism of the vector space such that $x\cdot V
\subset V$ then $x\cdot V = V$.

\end{cor}

%

Intersections of geometric hyperplanes are (special) linear
subspaces of $\Ext^k H$, and linear arrangements in $\Ext^k H$ are
finite unions of those linear subspaces. So, the previous results
apply, in particular, to the family of linear arrangements in the
exterior product. Let us call a linear arrangement $S$ in $\Ext^k
H$ \emph{economical} if it is contained in any other linear
arrangement $S'$ such that $\P S \cap\grass(k,H)=\P S'
\cap\grass(k,H)$. There is a natural bijection between economical
linear arrangements of the exterior product and linear
arrangements of the Grassmannian that preserves the inclusion
order. Hence, the previous results translate immediately to linear
arrangements in the Grassmannian. We summarize this in the next
corollary:

\begin{cor}\label{c.linear9} Both in $\Ext^k H$ and in $\grass(k,H)$,

\begin{enumerate}
\item The set of all linear arrangements is closed under finite unions
and arbitrary intersections.
\item Any totally ordered (under inclusion) set of linear arrangements is
well ordered.
\item If $S$ is a linear arrangement and $x$ is a linear isomorphism
such that $x \cdot S \subset S$ then $x \cdot S=S$.
\end{enumerate}

\end{cor}

\subsection{Symplectic spaces}

A {\it symplectic form} on a vector space $H$ is a bilinear form which is alternate, that
is,  $\omega(u,v)=-\omega(v,u)$ for all $u$ and $v$, and  non-degenerate, that is,
for every $u \in H \setminus \{0\}$ there exists $v \in H$ such that $\omega(u,v) \neq 0$.
We call $(H,\omega)$ a symplectic space.  Notice that $\dim H$ is necessarily an
even number $2g$.
A  {\it symplectic isomorphism} $A:(H,\omega) \to (H',\omega')$ is an isomorphism
satisfying $\omega(u,v)=\omega'(A \cdot u,A \cdot v)$.  By Darboux's theorem, all
symplectic spaces with the same dimension are symplectically isomorphic.

Given a subspace $F \subset H$, the {\it symplectic orthogonal} of $F$ is the set
$H_F$ of all $v \in H$ such that for every $u \in F$ we have $\omega(u,v)=0$.
Its dimension is complementary to the dimension of $F$, that is,
$\dim H_F=2g-\dim F$.  Notice that if two subspaces of $H$ intersect non-trivially
and have complementary dimension then their symplectic orthogonals also
intersect non-trivially.

A subspace $F \subset H$ is called {\it isotropic} if for every
$u,v \in F$ we have $\omega(u,v)=0$ or, in other words, if $F$ is
contained in its symplectic orthogonal. This implies that $\dim F
\leq g$. A {\it Lagrangian subspace} is an isotropic space of
maximal possible dimension $g$.  We represent by $\iso(k,H)
\subset \grass(k,H)$ the space of isotropic $k$-planes, for $1\le
k \le g$.  This is a closed, hence compact subset of the
Grassmannian. Notice that $\iso(1,H)=\grass(1,H)=\P H$. Given $F
\in \iso(k,H)$, we call  {\it symplectic reduction} of $H$ by $F$
the quotient space $H^F=H_F/F$.  Notice that $H^F$ admits a
canonical symplectic form $\omega^F$ defined by
$\omega^F([u],[v])=\omega(u,v)$.

Let $\LL(H)$ be the space of {\it Lagrangian flags}, that is, the
set of $(F_i)_{i=1}^g$ where $F_i$ is an isotropic subspace of $H$
of dimension $i$ and $F_i \subset F_{i+1}$ for all $i$.  There
exists a canonical embedding $\LL(H) \to \FF(H)$ given by
$(F_i)_{i=1}^g \mapsto (F_i)_{i=1}^{2g-1}$ where $F_i$ is the
symplectic orthogonal of $F_{2g-i}$. One can see $\LL(H)$ as a
fiber bundle over each $\iso(k,H)$, through the projection

 \be \label {fbs}
 \Upsilon_k:\LL(H) \to \iso(k,H),
 \quad (F_i)_{i=1}^g \mapsto F_k.
 \ee
We will use mostly the particular case $k=1$,
 \be \label {fiber bundle structure}
 \Upsilon=\Upsilon_1:\LL(H) \to \P H,
 \quad (F_i)_{i=1}^g  \mapsto F_1.
 \ee
The fiber over a given $\lambda\in\P H$ is naturally isomorphic to
$\LL(H^\lambda)$, via the isomorphism
 \be \label{fbsbis}
 \Upsilon^{-1}(\lambda) \to \LL(H^\lambda),
 \quad (F_i)_{i=1}^g \mapsto (F_{i+1}/\lambda)_{i=1}^{g-1}.
 \ee

\begin{lemma} \label {2.7}
Let $1 \leq k \leq g$.  Any hyperplane section in $\grass(k,H)$ intersects $\iso(k,H)$
in a non-empty compact subset with empty interior.
\end{lemma}

\begin{pf}
Let $S$ be a geometric hyperplane dual to some $E \in
\grass(2g-k,H)$. It is easy to see that $S$ meets $\iso(k,H)$:
just pick any $\lambda \in \P H$ contained in $E$ and consider any
element $F$ of $\iso(k,H)$ containing $\lambda$. By construction,
$F$ is in $S\cap\iso(k,H)$. Clearly, the intersection is closed in
$\iso(k,H)$, and so it is compact. We are left to show that the
complement of $S\cap \iso(k,H)$ is dense. If $k=1$ the result is
clear; in particular, this takes care of the case $g=1$. Next,
assuming the result is true for $(k-1,g-1)$ we deduce it is also
true for $(k,g)$. Indeed, for a dense subset $D$ of lines $\lambda
\in \P H$, we have that $\lambda \not \subset E$ and $E\not\subset
H_\lambda$. Consequently, $E_\lambda=(E\cap H_\lambda)/\lambda$ is
a subspace of $H^\lambda$ of dimension $2g-k-1$. By the induction
hypothesis, a dense subset $\DD_\lambda$ of $\iso(k-1,H^\lambda)$
is not in the hyperplane section $S_\lambda$ dual to $E_\lambda$.
Then, in view of \eqref{fbs}, the set
$\LL_\lambda=\Upsilon^{-1}_{k-1}(\DD_\lambda)$ of Lagrangian flags
whose $(k-1)$-dimensional subspace is contained in $\DD_\lambda$
is dense in $\LL(H^\lambda)$. By \eqref{fbsbis}, we may think of
$\LL_\lambda$ as a subset of $\LL(H)$ contained in the fiber over
$\lambda$. So, taking the union of the $\LL_\lambda$ over all
$\lambda\in D$ we obtain a dense subset $\LL$ of $\LL(H)$. Using
\eqref{fbs} once more, we conclude that $\Upsilon_k(\LL)$ is a
dense subset of $\iso(k,H)$. It suffices to prove that no element
of this set belongs to $S$. Let $F\in\Upsilon_k(\LL)$. By
definition, there exists $\lambda\subset F$ with
$\lambda\not\subset E$ and $(F/\lambda)\cap(E\cap
H_\lambda)/\lambda=\{0\}$. Since $F$ is contained in $H_\lambda$,
as it is isotropic, this implies that $F\cap E = \{0\}$, which is
precisely what we wanted to prove.
%
\end{pf}

\subsection{Linear actions of monoids}

A {\it monoid} $\BB$ is a set endowed with a binary operation that
is associative and admits a neutral element (the same axioms as in
the definition of group, except for the existence of inverse). The
monoids that interest us most in this context correspond to spaces
of loops through a given vertex in a Rauzy diagram, relative to
the concatenation operation. A {\it linear action} of a monoid is
an action by isomorphisms of a finite dimension vector space $H$.
It induces actions on the projective space $\P H$, the
Grassmannian spaces $\grass(k,H)$, and the flag space $\FF(H)$.
Given a subspace $F \subset H$, we denote by $\BB_F$ the
stabilizer of $F$, that is, the subset of $x\in\BB$ such that
$x\cdot F=F$.

A {\it symplectic action} of a monoid is an action by symplectic
isomorphisms on a symplectic space $(H,\omega)$. It induces
actions on $\iso(k,H)$ and the space of Lagrangian flags $\LL(H)$.
These actions are compatible with the fiber bundles \eqref {fbs},
in the sense that they are conjugated to each other by
$\Upsilon_k$. Given $\lambda \in \P H$, the stabilizer
$\BB_\lambda \subset \BB$ of $\lambda$ acts symplectically on
$H^\lambda$.

Let $\BB$ be a monoid acting by homeomorphisms of a compact space
$X$. By {\it minimal set} for the action of $\BB$ we mean a
non-empty closed set $C \subset X$ which is invariant, that is $x
\cdot C=C$ for all $x \in \BB$, and which has no proper subset
with those properties
\footnote {If $\BB$ is actually a group, a minimal set can be equivalently
defined as a non-empty closed set $C \subset X$ such that $\BB
\cdot x$ is dense in $C$ for every $x \in C$.}.  We say that the
action is {\it minimal} if the whole space is the only minimal
set. Detecting minimality of actions of a monoid in fiber bundles
can be reduced to detecting minimality for the action in the basis
and in the fiber.  We will need the following particular case of
this idea.

\begin{lemma} \label {2.2}
Let $\BB$ be a monoid acting symplectically on $(H,\omega)$.  Assume that
the action of $\BB$ on $\P H$ is minimal and that there exists $\lambda \in
\P H$ such that the stabilizer $\BB_\lambda \subset \BB$ of $\lambda$
acts minimally on $\LL(H^\lambda)$.  Then $\BB$ acts minimally on $\LL(H)$.
\end{lemma}

\begin{pf}
Let $C$ be a closed invariant subset of $\LL(H)$. For $\lambda'
\in \P H$, let $C_{\lambda'}$ be the intersection of $C$ with the
fiber over $\lambda'$. Notice that $C_\lambda$ may be seen as a
closed set invariant for the action of $\BB_\lambda$ on
$\LL(H^\lambda)$. So, $C_\lambda$ is either empty or the whole
fiber.  In the first case, let $\Lambda$ be the set of $\lambda'$
such that $C_{\lambda'}$ is non-empty. Then $\Lambda$ is a closed
invariant set which is not the whole $\P H$, and so it is empty.
This means $C$ itself is empty. In the second case, let $\Lambda$
be the set of $\lambda'$ such that $C_{\lambda'}$ coincides with
the whole fiber. This time, $\Lambda$ is a non-empty invariant set
and so it is the whole $\P H$. In other words, in this case
$C=\LL(H)$.
\end{pf}

\subsection{Singular values, Lyapunov exponents}

If we supply a vector space $H$ with an inner product, we identify $H$ with
the dual $H^*$, and we also introduce a metric on the Grassmannians: the
distance between $F,F' \in \grass(k,H)$ is the maximum of the angle between
lines $\lambda \subset F$ and $\lambda' \subset F'$.  We will often consider
balls with respect to this metric.  All balls will be assumed to have radius less
than $\pi/2$.

The inner product also allow us to speak of the {\it singular values} of a linear
isomorphism $x$ acting on $H$: those are the square roots of the eigenvalues
(counted with multiplicity) of the positive self-adjoint operator $x^*x$.
We always order them
\be
\sigma_1(x) \geq \cdots \geq \sigma_{\dim H}(x)>0.
\ee
A different inner product gives singular values differing from the $\sigma_i$ by
bounded factors, where the bound is independent of $x$.

The {\it Lyapunov exponents} of a linear isomorphism $x$ acting on $H$
are the logarithms of the absolute values of its eigenvalues, counted with
multiplicity. We denote and order the Lyapunov exponents of $x$ as
\be
\theta_1(x) \geq \cdots \geq \theta_{\dim H}(x).
\ee
Alternatively, they can be  defined by
\be
\theta_i(x)=\lim \frac {1} {n} \ln \sigma_i(x^n).
\ee

Given a linear isomorphism $x$ such that $\sigma_k(x)>\sigma_{k+1}(x)$, we
let $E^+_k(x)$ and $E^-_k(x)$ be the orthogonal spaces of dimension $k$ and
$\dim H-k$, respectively, such that $E^+_k(x)$ is spanned by the eigenvectors
of $x^*x$ with eigenvalue at least $\sigma_k(x)^2$ and $E^-_k(x)$ is spanned
by the eigenvectors of $x^*x$ with eigenvalue at most $\sigma_{k+1}(x)^2$.
For any $h=h^+_k+h^-_k \in E^+_k(x) \oplus E^-_k(x)$,
\be\label{eq.cotainferior}
\|x \cdot h\| \geq \sigma_k(x) \|h_k^+\|.
\ee
This is useful for the following reason.
If $F \in \grass(k,H)$ is transverse to $E^-_k(x)$ then
\be\label{projecao}
\|x \cdot h\| \geq c \sigma_k(x) \|h\| \quad\text{for all }  h \in F,
\ee
where $c>0$ does not depend on $x$ but only on the distance between
$F$ and the hyperplane section dual to $E^-_k(x)$ in $\grass(k,H)$.

We collect here several elementary facts from linear algebra that
will be useful in the sequel. Related ideas appear in Section 7.2
of \cite{AV}.

\begin{lemma} \label {compactlimit}
Let $x_n$ be a sequence of linear isomorphisms of $H$ such that
$\ln \sigma_k(x_n)-\ln \sigma_{k+1}(x_n) \to \infty$, and such
that $x_n \cdot E^+_k(x_n) \to E^u_k$ and $E^-_k(x_n) \to E^s_k$
as $n\to\infty$. If $K$ is a compact subset of $\grass(k,H)$ which
does not intersect the hyperplane dual to $E^s_k$ then $x_n \cdot
K \to E^u_k$ as $n\to\infty$.
\end{lemma}

\begin{pf}
Considering a subsequence if necessary, let $F_n \in K$ be such
that $x_n \cdot F_n$ converges to some $F' \neq E^u_k$. Take $h_n
\in F_n$ with $\|h_n\|=1$ such that
$$
\frac {x_n \cdot h_n}{\|x_n \cdot h_n\|} \to h \notin E^u_k.
$$
By \eqref{projecao}, there exists $c_1>0$ depending only on the
distance from $K$ to the hyperplane dual to $E_k^s$ such that
 \be
 \|x_n \cdot h_n\| \geq c_1 \sigma_k(x_n).
 \ee
Moreover, since $x_n \cdot E^+_k(x_n)=E^-_{\dim H-k}(x_n^{-1})$,
there exists $c_2>0$ depending only on the distance from $h$ to
$E_k^u$ such that for $n$ large
 \be
 1=\|h_n\|= \|x_n^{-1} x_n \cdot h_n\|
 \geq c_2 \sigma_{\dim H-k}(x_n^{-1}) \|x_n \cdot h_n\|.
 \ee
Consequently, since $\sigma_{\dim
H-k}(x_n^{-1})=\sigma_{k+1}(x_n)^{-1}$,
 \be
 1 \geq c_1c_2 \frac {\sigma_k(x_n)} {\sigma_{k+1}(x_n)}
 \ee
for all $n$, and this contradicts the hypothesis.
\end{pf}

\begin{lemma}\label {sigma theta}
Let $x_n$ be a sequence of linear isomorphisms of $H$ and let $1
\leq r \leq \dim H-1$ be such that $\ln \sigma_k(x_n)-\ln
\sigma_{k+1}(x_n) \to \infty$ for all $1 \leq k \leq r$. Assume
that $x_n \cdot E^+_k(x_n) \to E^u_k$, $E^-_k(x_n) \to E^s_k$, and
that $x$ is a linear isomorphism of $H$ such that $(x \cdot E^u_k)
\cap E^s_k=\{0\}$ for all $1 \leq k \leq r$. Then there exists
$C>0$ such that $|\theta_k(x x_n)-\ln \sigma_k(x_n)|<C$ for all $1
\leq k \leq r$ and $n$ large.
\end{lemma}

\begin{pf}
Let $U_k$ be an open ball around $x \cdot E^u_k$ such that $F \cap E^s_k=\{0\}$
for every $F \in \overline U_k$.  By the previous lemma, for $n$ large we have
$x x_n \cdot \overline U_k \subset U_k$.  In particular, there exists
$E^u_{k,n} \in U_k$ such that $x x_n \cdot E^u_{k,n}=E^u_{k,n}$ and
 \be
\|x x_n \cdot h\| \geq c_1 \sigma_k(x_n) \|h\| \quad\text{for all
} h \in E^u_{k,n},
 \ee
where $c_1$ depends only on $x$ and the distance from $\overline
U_k$ to $E_k^s$. Consequently,
 \be \label{eq.2.17}
 e^{\theta_k(x x_n)} \geq c_1 \sigma_k(x_n)
 \quad \text{for all } 1 \leq k \leq r.
 \ee
Clearly, we also have
 \be
 \prod_{1 \leq k \leq j} e^{\theta_k(x x_n)}
 \leq \prod_{1 \leq k \leq j} \sigma_k(x x_n)
 \leq c_2^{-1} \prod_{1 \leq k \leq j} \sigma_k(x_n)
 \quad\text{for all } 1 \leq j \leq r,
 \ee
where $c_2=c_2(x)$. Using \eqref{eq.2.17}, we conclude that
 \be
 e^{\theta_j(x x_n)} \leq c^{-j} \sigma_j(x_n)
 \quad\text{for all } 1 \leq j \leq r,
 \ee
with $c=\min\{c_1,c_2\}$, and the result follows.
\end{pf}

\begin{lemma} \label{l.2.11}
Let $x_n$ be a sequence of linear isomorphisms of $H$. Suppose there is
$F \in \grass(k,H)$ such that the set $\{F' \in \grass(k,H),\, x_n \cdot F' \to F\}$ is
not contained in a hyperplane section.
Then we have $\ln \sigma_k(x_n)-\ln \sigma_{k+1}(x_n) \to \infty$ and
$F=\lim x_n \cdot E^+_k(x_n)$.
\end{lemma}

\begin{pf}
Assume that $\ln \sigma_k(x_n)-\ln \sigma_{k+1}(x_n)$ is bounded
(along some subsequence). Passing to a subsequence, and replacing
$x_n$ by $y_n x_n z_n$, with $y_n, y_n^{-1}, z_n, z_n^{-1}$
bounded, we may assume that there exists $l \leq k<r$ such that
$\sigma_j(x_n)=\sigma_k(x_n)$ if $l \leq j \leq r$ and $|\ln
\sigma_j(x_n)-\ln \sigma_k(x_n)| \to \infty$ otherwise, and there
exists an orthonormal basis $\{e_i\}_{i=1}^{\dim H}$, independent
of $n$, such that $x_n \cdot e_i=\sigma_i(x_n) e_i$ for every $i$
and $n$. Let $E^u$, $E^c$, and $E^s$ be the spans of
$\{e_i\}_{i=1}^{l-1}$, $\{e_i\}_{i=l}^r$, and
$\{e_i\}_{i=r+1}^{\dim H}$ respectively. Notice that, for any
$F'\in \grass(k,H)$
\begin{enumerate}
\item $(E^c \oplus E^s)+F'=H$ if and only if any limit of $x_n \cdot F'$
      contains $E^u$,
\item $F' \cap E^s=\{0\}$ if and only if any limit of $x_n \cdot F'$
      is contained in $E^u \oplus E^c$,
\item If both (1) and (2) hold then $x_n \cdot F' \to F$ where $F=E^u \oplus
F^c$ and $F^c \subset E^c$ is such that $F^c \oplus E^s=F' \cap (E^c \oplus
E^s)$.
\end{enumerate}
Thus, if $x_n \cdot F' \to F$ then $F'$ is contained in  the
hyperplane section dual to $G \in \grass(\dim H-k,H)$ chosen as
follows.  If $E^u \not \subset F$ then $G$ can be any subspace
contained in $E^c \oplus E^s$. If $F \not \subset E^u \oplus E^c$
then $G$ can be any subspace containing $E^s$. If $E^u \subset F
\subset E^c \oplus E^u$ then $G$ can be any subspace containing
$E^s$ and such that $G \cap F \cap (E^c \oplus E^s) \neq \{0\}$.
This shows that if $\{F' \in \grass(k,H),\, x_n \cdot F' \to F\}$
is not contained in a hyperplane section then the difference $\ln
\sigma_k(x_n)-\ln\sigma_{k+1}(x_n) \to \infty$. To conclude, we
may assume that $E^-_k(x_n)$ converges to some $E^s_k$.  By
Lemma~\ref {compactlimit}, if $\{F' \in \grass(k,H),\, x_n \cdot
F' \to F\}$ is not contained in the hyperplane section dual to
$E^s_k$ then $F=\lim x_n \cdot E^+(x_n)$.
\end{pf}

\begin{lemma} \label {non-trivial}
Let $x_n$ be a sequence of linear isomorphisms of $H$, $\rho$ be a probability measure
on $\grass(k,H)$, and $\rho^{(n)}$ be the push-forwards of $\rho$ under $x_n$.
\begin{enumerate}
\item Assume that $\rho$ gives zero weight to any hyperplane section.  If
$\ln \sigma_k(x_n)-\ln\sigma_{k+1}(x_n) \to \infty$ and $x_n \cdot
E^+_k(x_n) \to E^u_k$, then $\rho^{(n)}$ converges in the weak$^*$
topology to a Dirac mass on $E^u_k$.
\item Assume that $\rho$ is not supported in a hyperplane section. If $\rho^{(n)}$
converges in the weak$^*$ topology to a Dirac mass on $E^u_k$ then
$\ln \sigma_k(x_n)-\ln\sigma_{k+1}(x_n) \to \infty$ and $x_n \cdot E^+_k(x_n) \to E^u_k$.
\end{enumerate}
\end{lemma}

\begin{pf}
Let us first prove (1). We may assume that $E^-_k(x_n)$ also
converges to some $E^s_k$. Take a compact set $K$ disjoint from
the hyperplane section dual to $E^s_k$ and such that
$\rho(K)>1-\epsilon$. Then $\rho^{(n)}(x_n \cdot K)>1-\epsilon$
and $x_n \cdot K$ is close to $E^u_k$ for all large $n$, by
Lemma~\ref{compactlimit}. This shows that $\rho^{(n)}$ converges
to the Dirac measure on $E_k^u$. Now let us prove (2). The
hypothesis implies that, passing to a subsequence of an arbitrary
subsequence, $x_n \cdot F \to E^u_k$ for a full measure set, and
this set is not contained in a hyperplane section. Using
Lemma~\ref{l.2.11}, we conclude that $\ln \sigma_k(x_n)-\ln
\sigma_{k+1}(x_n) \to \infty$ and $x_n \cdot E^+_k(x_n) \to
E^u_k$.
\end{pf}

\begin{lemma} \label {fthetaargument}
Let $F^u \in \grass(k,H)$ and $F^s \in \grass(\dim H-k,H)$ be orthogonal subspaces.
Assume that $x$ is a linear isomorphism of $H$ such that $x \cdot F^u=F^u$ and
$x \cdot F^s=F^s$ and there is an open ball $U$ around $F^u$ such that
$x \cdot \overline U \subset U$.  Then $\sigma_k(x)>\sigma_{k+1}(x)$ and
$E^+_k(x)=F^u$ and $E^-_k(x)=F^s$.
\end{lemma}

\begin{pf}
Up to composition with orthogonal operators preserving $F^u$ and $F^s$, we
may assume that $x$ is diagonal with respect to some orthonormal basis
$e_1, \ldots ,e_{\dim H}$, and that $F^u$ and $F^s$ are spanned by elements of
the basis.  We may order the elements of the basis so that $x \cdot e_i=
\sigma_i(x) e_i$.  Let $F^u$ be spanned by $\{e_{i_j}\}_{j=1}^k$.
If $\sigma_l(x) \leq \sigma_r(x)$ for some $l \in \{i_j\}_{j=1}^k$,
$r \notin \{i_j\}_{j=1}^k$, let $F_\theta$ be spanned by
the $e_s$ which are either of the form $e_{i_j}$ with $i_j \neq l$ or of the
form $\cos 2 \pi \theta e_l+\sin 2 \pi \theta e_r$.  If
$\sigma_l(x)=\sigma_r(x)$ then $x \cdot F_\theta=F_\theta$ for every
$\theta$ and if $\sigma_l(x)<\sigma_r(x)$ then $x^n \cdot F_\theta \to
F_{\pi/2} \notin U$ for every $\theta$ such that $F_\theta \neq F_0$.
In both cases, this gives a contradiction (by considering some $F_\theta \in
\partial U$).
\end{pf}

If $(H,\omega)$ is a symplectic space (with $\dim H=2g$), the choice of an
inner product defines an antisymmetric linear isomorphism $\Omega$ on $H$
satisfying
\be
\omega(\Omega \cdot u,\Omega \cdot v)=\langle u,\Omega \cdot v \rangle.
\ee
One can always choose the inner product so that $\Omega$ is also orthogonal
(we call such an inner product {\it adapted} to $\omega$).
If a linear isomorphism $x$ is symplectic, we have, for this particular
choice of inner product,
\be
\sigma_i(x) \sigma_{2g-i+1}(x)=1 \quad\text{for all } i=1, \ldots ,g
\ee
and hence
\be
\theta_i(x)=-\theta_{2g-i+1}(x)  \quad\text{for all } i=1, \ldots ,g.
\ee

\begin{lemma} \label {isolimit}
Let $x_n$ be a sequence of symplectic isomorphisms of $H$ such that
$\sigma_k(x_n) \to \infty$ and $\sigma_k(x_n)>\sigma_{k+1}(x_n)$ and such
that $x_n \cdot E^+_k(x_n)$ converges to some space $E^u_k$.  Then $E^u_k$
is isotropic.
\end{lemma}

\begin{pf}
By \eqref{eq.cotainferior} , if $h_n \in x_n \cdot E^+_k(x_n)$ then
$\|x_n^{-1} \cdot h_n\| \leq \sigma_k(x_n)^{-1} \|h_n\|$.
Thus, if $u_n, v_n \in x_n \cdot E^+_k(x_n)$ are such that $\|u_n\|=\|v_n\|=1$, then
\be
|\omega(u_n,v_n)| = |\omega(x_n^{-1} \cdot u_n, x_n^{-1} \cdot v_n)|
 \leq c^{-1} \sigma_k(x_n)^{-2}.
\ee
Passing to the limit as $n\to\infty$, we get $\omega(u,v)=0$ for every
$u, v \in E^u_k$.
\end{pf}

\section{Rauzy classes and the Zorich cocycle} \label {rauzy}

\subsection{Interval exchange transformations}

We follow the presentation of \cite {MMY}.  An interval exchange
transformation is defined as follows. Let $\AA$ be some fixed
alphabet on $d\ge 2$ symbols. All intervals will be assumed to be
closed on the left and open on the right.
\begin{itemize}
\item Take an interval $I\subset\R$ and break it into subintervals
$\{I_x\}_{x \in \AA}$.
\item Rearrange the intervals in a new
order, via translations, inside $I$.
\end{itemize}

Modulo translations, we may always assume that the left endpoint of $I$ is $0$.
Thus the interval exchange transformation is entirely defined by the following data:
\begin{enumerate}
\item The lengths of the intervals $\{I_x\}_{x \in \AA}$,
\item Their orders before and after rearranging.
\end{enumerate}
The first is called length data, and is given by a vector $\lambda \in
\R^\AA_+$.  The second is called combinatorial data, and is given by a pair
of bijections $\pi=(\pi_0,\pi_1)$ from $\AA$ to $\{1,\ldots,d\}$ (we will
sometimes call such a pair of bijections a permutation).
We denote the set of all such pairs of bijections by $\sssigma(\AA)$.
We view a bijection $\AA \to \{1,\ldots,d\}$ as a row where the elements of $\AA$ are
displayed in the right order.  Thus we can see an element of $\sssigma(\AA)$ as
a pair of rows, the top (corresponding to $\pi_0$) and the bottom (corresponding to
$\pi_1$) of $\pi$.  The interval exchange transformation associated to this data will
be denoted $f=f(\lambda,\pi)$.

Notice that if the combinatorial data is such that the set of the
first $k$ elements in the top and bottom of $\pi$ coincide for
some $1 \leq k<d$ then, irrespective of the length data, the
interval exchange transformation  splits into two simpler
transformations. Thus we will consider only combinatorial data for
which this does not happen, which we will call ${\it
irreducible}$. Let $\ssigma(\AA) \subset \sssigma(\AA)$ be the set
of irreducible combinatorial data.

\subsection{Translation vector}

The positions of the intervals $I_x$ before and after applying the interval
exchange transformation differ by translations by
\be
\delta_x=\sum_{\pi_1(y)<\pi_1(x)} \lambda_y-\sum_{\pi_0(y)<\pi_0(x)}
\lambda_y.
\ee
We let $\delta=\delta(\lambda,\pi) \in \R^\AA$ be the {\it translation vector},
whose coordinates are given by the $\delta_x$.
Notice that the ``average translation''
$\langle \lambda,\delta \rangle = \sum_{x\in^\AA} \lambda_x \delta_x$
is zero. We can write the relation between $\delta$ and $(\lambda,\pi)$ as
\be
\delta(\lambda,\pi)=\Omega(\pi) \cdot \lambda,
\ee
where $\Omega(\pi)$ is a linear operator on $\R^\AA$ given by
\be\label{defOmega}
\langle \Omega(\pi) \cdot e_x,e_y \rangle=
\left \{ \begin{array}{ll}
1, & \pi_0(x)>\pi_0(y), \pi_1(x)<\pi_1(y),\\[5pt]
-1, & \pi_0(x)<\pi_0(y), \pi_1(x)>\pi_1(y),\\[5pt]
0, & \text {otherwise}.
\end{array}
\right .
\ee
(here $e_x$ is the canonical basis of $\R^\AA$ and the inner product
$\langle\, ,\, \rangle$ is the natural one which makes this canonical basis
orthonormal).
Notice that $\Omega(\pi)$ can be viewed as an antisymmetric matrix with
integer entries.  Notice also that $\Omega(\pi)$ may not be invertible.

We denote $H(\pi)=\Omega(\pi) \cdot \R^\AA$, which is the space
spanned by all possible translation vectors $\delta(\lambda,\pi)$.
We define a symplectic form $\omega=\omega_\pi$ on $H(\pi)$ by
putting \be \omega_\pi(\Omega(\pi) \cdot u,\Omega(\pi) \cdot
v)=\langle u,\Omega(\pi) \cdot v \rangle. \ee We let $2g(\pi)$ be
the dimension of $H(\pi)$, where $g(\pi)$ is called the {\it
genus}.

\subsection{Rauzy diagrams and monoids}

A {\it diagram} (or directed graph) consists of two kinds of
objects, vertices and (oriented) arrows joining two vertices.
Thus, an arrow has a start and an end. A {\it path} in the diagram
of length $m \geq 0$ is a sequence $v_0, \ldots, v_m$ of vertices
and a sequence of arrows $a_1, \ldots, a_m$ such that $a_i$ starts
at $v_{i-1}$ and ends in $v_i$.  If $\gamma_1$ and $\gamma_2$ are
paths such that the end of $\gamma_1$ is the start of $\gamma_2$,
their concatenation is also a path, denoted by $\gamma_1\gamma_2$.
The set of all paths starting and ending at a given vertex $v$ is
a monoid for the operation of concatenation. We can identify paths
of length zero with vertices and paths of length one with arrows.

Given $\pi \in \ssigma(\AA)$ we consider two operations.  Let $x$
and $y$ be the last elements of the top and bottom rows. The  {\it
top}  operation keeps the top row; on the other hand, it takes $y$
and inserts it back in the bottom immediately to the right of the
position occupied by $x$. When applying this operation to $\pi$,
we will say that $x$  {\it  wins} and $y$ {\it loses}. The  {\it
bottom} operation is defined in a dual way, just interchanging the
words top and bottom, and the roles of $x$ and $y$. In this case
we say that $y$ wins and  $x$ loses. Notice that both operations
preserve the first elements of both the top and the bottom row.

It is easy to see that those operations give bijections of
$\ssigma(\AA)$. The  {\it Rauzy diagram} associated to $\AA$ has
the elements of $\ssigma(\AA)$ as its vertices, and its arrows
join each vertex to the ones obtained from it through either of
the operations we just described. So, every vertex is the start
and end of two arrows, one top and one bottom. Thus, every arrow
has a start, an end, a type (top/bottom), a winner and a loser.
The set of all paths is denoted by $\Pi(\AA)$.

The orbit of $\pi$ under the monoid generated by the actions of
the top and bottom operations on $\ssigma(\AA)$ will be called the
{\it Rauzy class} of $\pi$, and denoted $\RR(\pi)$.  The set of
all paths inside a given Rauzy class will be denoted $\Pi(\RR)$.
The set of all paths that begin and end at $\pi\in \ssigma(\AA)$
will be denoted $\Pi(\pi)$.  We call $\Pi(\pi)$ a {\it Rauzy
monoid}.

\subsection{Rauzy induction}

Let $\RR \subset \ssigma(\AA)$ be a Rauzy class, and define
$\Delta^0_\RR=\R^\AA_+ \times \RR$. Given $(\lambda,\pi)$ in
$\Delta^0_\RR$, we say that we can apply Rauzy induction to
$(\lambda,\pi)$ if $\lambda_x \neq \lambda_y$, where $x, y \in
\AA$ are the last elements of the top and bottom rows of $\pi$,
respectively. Then we define $(\lambda',\pi')$ as follows:
\begin{enumerate}
\item Let $\gamma=\gamma(\lambda,\pi)$ be a top or bottom arrow on the Rauzy
diagram starting at $\pi$, according to whether $\lambda_x>\lambda_y$ or $\lambda_y>\lambda_x$.
\item Let $\lambda'_z=\lambda_z$ if $z$ is not the winner of $\gamma$, and
$\lambda_z=\max \{\lambda_x,\lambda_y\}-\min \{\lambda_x,\lambda_y\}$ if $z$
is the winner of $\gamma$,
\item Let $\pi'$ be the end of $\gamma$.
\end{enumerate}
We say that $(\lambda',\pi')$ is obtained from $(\lambda,\pi)$ by
applying Rauzy induction, of type top or bottom depending on
whether the type of $\gamma$ is top or bottom. We have that $\pi'
\in \ssigma(\AA)$ and $\lambda' \in \R^\AA_+$. The interval
exchange transformations $f:I \to I$ and $f':I' \to I'$ specified
by the data $(\lambda,\pi)$ and $(\lambda',\pi')$ are related as
follows.  The map $f'$ is the first return map of $f$ to a
subinterval of $I$, obtained by cutting from $I$ a subinterval
with the same right endpoint and of length $\lambda_z$, where $z$
is the loser of $\gamma$. The map $Q_R:(\lambda,\pi) \to
(\lambda',\pi')$ is called  {\it Rauzy induction map}. Its domain
of the definition, the set of all $(\lambda,\pi) \in \Delta^0_\RR$
such that $\lambda_x \neq \lambda_y$, will be denoted by
$\Delta^1_\RR$.

\subsection{Relation between translation vectors}

Let us calculate the relation between the translation vectors $\delta=\delta(\lambda,\pi)$
and $\delta'=\delta(\lambda',\pi')$.
Let $\gamma$ be the arrow that starts at $\pi$ and ends at $\pi'$. We have
$\delta'_z=\delta_z$ if $z$ is not the loser of $\gamma$ and $\delta'_z=\delta'_x+\delta'_y$
if $z$ is the loser of $\gamma$ (we continue to denote by $x, y$ the last symbols in the
top and bottom rows of $\pi$, respectively).  In other words, we can write
\be
\delta(\lambda',\pi')=\Theta(\gamma) \cdot \delta(\lambda,\pi),
\ee
where $\Theta(\gamma)$ is the linear operator of $\R^\AA$ defined by
$\Theta(\gamma) \cdot e_z=e_x+e_y$ if $z$ is the winner of $\gamma$,
and $\Theta(\gamma) \cdot e_z=e_z$ otherwise. Notice that we also have
\be
\lambda=\Theta(\gamma)^* \cdot \lambda',
\ee
where the adjoint $\Theta(\gamma)^*$ is taken with respect to the natural
inner product on $\R^\AA$ that renders the canonical basis orthonormal.
Consequently,  if $h \in \R^\AA$ and $h'=\Theta(\gamma) \cdot h$ then
$\langle \lambda,h \rangle=\langle \lambda',h' \rangle$.

Since $H(\pi)$ and $H(\pi')$ are spanned by possible translation
vectors, $\Theta(\gamma) \cdot H(\pi)=H(\pi')$, and so the
dimension of $H(\pi)$ only depends on the Rauzy class of $\pi$.
One can check that \be\label{eq.omegasymplectic} \Theta(\gamma)
\Omega(\pi) \Theta(\gamma)^*=\Omega(\pi'), \ee which implies that
$\Theta(\gamma):(H(\pi),\omega_\pi) \to (H(\pi'),\omega_{\pi'})$
is a symplectic isomorphism. Notice that $\Theta(\gamma)$ can be
viewed as a matrix with non-negative integer entries and
determinant $1$.  We extend the definition of $\Theta$ from arrows
to paths in $\Pi(\AA)$ in the natural way, $\Theta(\gamma_1
\gamma_2)=\Theta(\gamma_2) \Theta(\gamma_1)$. In this way $\Theta$
induces a representation on $\SL(\AA,\Z)$ of the Rauzy monoids
$\Pi(\pi) \subset \Pi(\AA)$, $\pi \in \ssigma(\AA)$.

\subsection{Iterates of Rauzy induction}

The connected components of $\Delta^0_\RR=\R_+^\AA\times\RR$ are
naturally labelled by the elements of $\RR$ or, in other words, by
length $0$ paths in $\Pi(\RR)$ .  The connected components of the
domain $\Delta^1_\RR$ of the induction map $Q_R$ are naturally
labelled by arrows, that is, length $1$ paths in $\Pi(\RR)$.  One
easily checks that each connected component of $\Delta^1_\RR$ is
mapped bijectively to some connected component of $\Delta^0_\RR$.
Now let $\Delta^n_\RR$ be the domain of $Q_R^n$, for each $n \geq
2$. The connected components of $\Delta^n_\RR$ are naturally
labelled by length $n$ paths in $\Pi(\RR)$: if $\gamma$ is
obtained by concatenation of arrows $\gamma_1, \ldots, \gamma_n$,
then $\Delta_\gamma=\{x \in \Delta^0(\lambda,\pi) : Q_R^{k-1}(x)
\in \Delta_{\gamma_k},\, 1 \leq k \leq n\}$. If $\gamma$ is a
length $n$ path in $\Pi(\RR)$ ending at $\pi \in \RR$, then \be
Q^\gamma_R=Q^n_R:\Delta_\gamma \to \Delta_\pi \ee is a bijection.

The set $\Delta_\RR=\cap_{n \geq 0} \Delta^n_\RR$ of
$(\lambda,\pi)$ to which we can apply Rauzy induction infinitely
many times\footnote{There is a nice explicit characterization of
this set, called the Keane property. See \cite {MMY} for a
statement.} contains all $(\lambda,\pi)$ such that the coordinates
of $\lambda$ are rationally independent, and so it is a full
Lebesgue measure subset of $\Delta^0_\RR$. Let $(\lambda,\pi) \in
\Delta_\RR$.  The connected components of $\Delta^m_\RR$, $m \geq
0$ are a nested sequence of convex cones containing the half-line
$\{(t\lambda,\pi),\, t \in \R_+\}$, and their intersection is the
connected component of $(\lambda,\pi)$ in $\Delta_\RR$. In
general, it is not true that this connected component reduces to
the half-line $\{(t\lambda,\pi),\, t \in \R_+\}$ (``combinatorial
rigidity''): this happens precisely when the interval exchange
transformation  defined by $(\lambda,\pi)$ is uniquely ergodic,
and one can find counterexamples as soon as the genus $g=\dim
H(\pi)/2$ is at least $2$.

\subsection{Rauzy renormalization map}

Since $Q_R$ commutes with dilations, it projectivizes to a map
$R_R:\P\Delta^1_\RR \to \P\Delta^0_\RR$, that we call {\it  Rauzy renormalization
map}.  Let  $R^\gamma_R$ be the projectivization of $Q^\gamma_R$, for
each path $\gamma$.

\begin{thm}[Masur \cite {M}, Veech \cite {V1}] \label {MV}
The map $R_R$ has an ergodic conservative absolutely continuous invariant infinite
measure $\mu_R$. The density is analytic and positive in $\P\Delta^0_\RR$.
\end{thm}

Conservativeness means that if a measurable set contains its pre-image then the
difference between the two has zero measure. It ensures that Poincar\'e recurrence
holds in this context, despite the fact that the measure is infinite.
This theorem is the key step in the proof by Masur and Veech that almost every interval
exchange transformation  is uniquely ergodic. Indeed, they show that if $x$ is recurrent
under $R_R$, then the connected component of $x$ in $\P\Delta_\RR$ reduces to a
point. The latter implies unique ergodicity of the interval exchange transformation  given
by the $(\lambda,\pi)$ that projectivize to $x$.

\subsection{The Zorich map}

The Rauzy renormalization map does not admit an absolutely
continuous invariant probability because it is too slow. For
instance, in the case of two intervals, the Rauzy renormalization
map is just the Farey map, which exhibits a parabolic fixed point.
Zorich introduced a way to ``accelerate'' Rauzy induction, that
produces a new renormalization map, which is more expanding, and
always admits an absolutely continuous invariant probability. In
the case of two intervals, the Zorich renormalization map is,
essentially, the Gauss map.

We say that we can apply Zorich induction to some $(\lambda,\pi)
\in \Delta^0_\RR$ if there exists a smallest $m \geq 1$ such that
we can apply Rauzy induction $m+1$ times to $(\lambda,\pi)$, and
in doing so we use both kinds of operations, top and bottom.  Then
we define $Q_Z(\lambda,\pi)=Q_R^m(\lambda,\pi)$. The domain of
this {\it Zorich induction map} is the union of $\Delta_\gamma$
over all paths $\gamma$ of length $m+1\geq 2$ which are obtained
by concatenating $m$ arrows of one type (top or bottom) followed
by an arrow of the other type. If we can apply Rauzy induction
infinitely many times to $(\lambda,\pi)$ then we can also apply
Zorich induction infinitely many times. The projectivization of
the Zorich induction map $Q_Z$ is called the {\it Zorich
renormalization map} $R_Z$ or, simply, Zorich map.

\begin{thm}[Zorich \cite {Z2}] \label {Z}
The Zorich map $R_Z$ admits an ergodic absolutely continuous
invariant probability measure. The density is analytic and
positive in $\P\Delta^1_\RR$.
\end{thm}

We call the measure given by this theorem the {\it Zorich measure} $\mu_Z$.

\subsection{The Zorich cocycle}

We define a linear cocycle $(x,h) \mapsto (R_Z(x),B^Z(x) \cdot h)$
over the Zorich map $x \mapsto R_Z(x)$, as follows. Let $x$ belong
to a connected component $\P\Delta_\gamma$ of the domain of $R_Z$,
where $\gamma$ is a length $m+1\geq 2$ path, and let $\tilde
\gamma$ be the length $m$ path obtained by dropping the very last
arrow in $\gamma$ (the one that has type distinct from all the
others). Then $B^Z(x)=\Theta(\tilde \gamma)$. To specify the
linear cocycle completely we also have to specify where $h$ is
allowed to vary.  There are two natural possibilities: either
$H(\pi)$ for $x \in \P\Delta_\pi$, or the whole $\R^\AA$. In the
first case we speak of the {\it Zorich cocycle}, whereas in the
second one we call this the {\it extended Zorich cocycle}.

\begin{thm}[Zorich \cite {Z2}]
The (extended) Zorich cocycle is measurable:
\be
\int \ln \|B^Z(x)\| d\mu_Z(x)<\infty.
\ee
\end{thm}

The relation \eqref{eq.omegasymplectic} gives that the Zorich
cocycle is symplectic. Consequently, its Lyapunov exponents
$\theta_1 \geq \cdots \geq \theta_{2g}$ satisfy
$\theta_i=-\theta_{2g-i+1}$ for all $i$, where $g=g(\RR)$ is the
genus. We say that the Lyapunov spectrum is symmetric.  The
Lyapunov spectrum of the extended Zorich cocycle consists of the
Lyapunov spectrum of the Zorich cocycle together with additional
zeros.

\subsection{Inverse limit}

Given $\pi=(\pi_0,\pi_1)$, let $\Gamma_\pi \subset \R^\AA$ be the
set of all $\tau$ such that
 \be
 \sum_{\pi_0(z) \leq k} \tau_z>0 \quad\text{and}\quad
 \sum_{\pi_1(z) \leq k} \tau_z<0 \quad\text{for all} \quad
 1 \leq k \leq d-1.
 \ee
Notice that $\Gamma_\pi$ is an open cone. If $\gamma$ is the top
arrow ending at $\pi'$, let $\Gamma_\gamma$ be the set of all
$\tau \in \Gamma_{\pi'}$ such that $\sum_{x \in \AA} \tau_x<0$,
and if $\gamma$ is a bottom arrow ending at $\pi'$, let
$\Gamma_\gamma$ be the set of all $\tau \in \Gamma_{\pi'}$ such
that $\sum_{x \in \AA} \tau_x>0$. If $\gamma$ is an arrow starting
at $\pi$ and ending at $\pi'$ then \be \Theta(\gamma)^* \cdot
\Gamma_\gamma=\Gamma_\pi. \ee Thus, the map
$$
\widehat Q^\gamma_R:
 \Delta_\gamma \times \Gamma_\pi \to \Delta_{\pi'} \times \Gamma_\gamma, \quad
 \widehat Q^\gamma_R(x,\tau)=(Q_R(x),(\Theta(\gamma)^{-1})^* \cdot \tau)
$$
is invertible. Now we can define an invertible skew-product $\widehat Q_R$ over
$Q_R$ by putting together the $\widehat Q^\gamma_R$ for every arrow $\gamma$.
This is a map from $\cup \Delta_\gamma \times \Gamma_\pi$ (where the union
is taken over all $\pi \in \RR$ and all arrows $\gamma$ starting at $\pi$) to
$\cup \Delta_{\pi'} \times \Gamma_\gamma$ (where the union is taken over
all $\pi' \in \RR$ and all arrows ending at $\pi'$).

Let $\widehat \Delta_\RR^0=\cup_{\pi \in \RR} \Delta_\pi \times
\Gamma_\pi$, and let $\widehat \Delta_\RR \subset \widehat
\Delta_\RR^0$ be the set of all points that can be iterated
infinitely many times forward and backwards by $\widehat Q_R$.
Note that $((\lambda,\pi),\tau) \in \widehat \Delta^0_\RR$ can be
iterated infinitely many times forward/backwards by $\widehat Q_R$
provided that the coordinates of $\lambda/\tau$ are rationally
independent. Projectivization in $\lambda$ and $\tau$ gives a map
$\widehat R_R:\widehat {\P\Delta}_\RR \to \widehat
{\P\Delta}_\RR$. This is an invertible map that can be seen as the
inverse limit of $R_R$: the connected components of $\widehat
{\P\Delta}_\RR$ reduce to points Lebesgue almost everywhere.
There is a natural infinite ergodic conservative invariant measure
$\hat \mu_R$ for $\widehat R_R$ that is equivalent to Lebesgue
measure. The projection $\mu_R$ of $\hat \mu_R$ is the measure
appearing in Theorem \ref {MV}: this is how $\mu_R$ is actually
constructed in \cite {V1}.

Let $\Upsilon_\RR^0=\cup \P\Delta_\gamma \times
\P\Gamma_{\gamma'}$, where the union is taken over all pairs of
arrows $\gamma'$ and $\gamma$ of distinct types (one is top and
the other is bottom) such that $\gamma'$ ends at the start of
$\gamma$. 
Let $\Upsilon_\RR=\Upsilon_\RR^0 \cap \widehat{\P\Delta}_\RR$ and
$\widehat R_Z:\Upsilon_\RR \to \Upsilon_\RR$ be the first return
map. The choice of $\Upsilon_\RR$ is so that $\widehat
\mu_R(\Upsilon_\RR)$ is finite.  So, we may normalize
$\hat\mu_R|\Upsilon_\RR$ to get a probability measure $\hat \mu_Z$
on $\Upsilon_\RR$ which is invariant under $\widehat R_Z$.  One
checks that $\hat \mu_Z$ is ergodic. Moreover,  the map $\widehat
R_Z$ is a skew-product over $R_Z$ that can be seen as the inverse
limit of $R_R$.  The projection $\mu_Z$ of $\hat \mu_Z$ is the
probability described in Theorem \ref {Z}: this is how $\mu_Z$ is
constructed in \cite {Z2}.

We call {\it invertible Zorich cocycle} the lift $(\widehat R_Z,\widehat B^Z)$ of the
Zorich cocycle to a cocycle over the invertible Zorich map, defined by
$\widehat B^Z((\lambda,\pi),\tau)=B^Z(\lambda,\pi)$.

\subsection{Minimality in the projective space}

The aim of this section is to prove Corollary~\ref {k=1}, regarding the projective
behavior of the matrices involved in the Zorich cocycle.  This result follows easily
from known constructions, but will be quite useful in the sequel .

Given $x \in \widehat{\P\Delta}_\RR$, obtained by projectivizing
$((\lambda,\pi),\tau)$, consider the subspaces $E^{uu}$, $E^c$,
$E^{ss}$ of $H(\pi)$ defined as follows: firstly, $E^{uu}$ is the
line spanned by $\Omega(\pi) \cdot \tau$; secondly, $E^{ss}$ is
the line spanned by $\Omega(\pi) \cdot \lambda$; and, finally,
$E^c$ is the symplectic orthogonal to the plane $E^{uu} \oplus
E^{ss}$. Notice that $E^{uu}$ is not symplectically orthogonal to
$E^{ss}$: $\langle \lambda,\Omega(\pi) \cdot \tau \rangle<0$,
since both $\lambda$ and $-\Omega(\pi) \cdot \tau$ have only
positive coordinates. Thus, $E^c$ has codimension $2$ in $H(\pi)$
and $H(\pi)=E^{uu} \oplus E^c \oplus E^{ss}$.

\begin{lemma}
The splitting $E^{uu} \oplus E^c \oplus E^{ss}$ is invariant under the invertible Zorich cocycle.
The spaces $E^{uu}$ and $E^{ss}$ correspond to the largest and the smallest Lyapunov
exponents, respectively, and $E^c$ corresponds to the remaining exponents.
\end{lemma}

\begin{pf}
The invariance of $E^{uu}$ and $E^{ss}$ follows directly from the
definitions and \eqref{eq.omegasymplectic}. The invariance of
$E^c$ is a consequence, since the Zorich cocycle is symplectic.
Since $E^{uu}$ is the projectivization of a direction in the
positive cone and the matrices of the invertible Zorich cocycle
have non-negative entries, $E^{uu}$ must be contained in the
subspace corresponding to the largest Lyapunov exponent.  Since
the largest Lyapunov exponent is simple (see \cite {Z2}), $E^{uu}$
must be the Oseledets direction associated to it. This implies
that $E^{ss}$ corresponds to the smallest Lyapunov exponent, since
$E^{ss}$ is an invariant direction which is not contained in the
symplectic orthogonal to $E^{uu}$.  It follows that $E^c$ must
correspond to the other Lyapunov exponents.
\end{pf}

It follows from the definitions that $-\Omega(\pi) \cdot
\Gamma_\pi\subset H(\pi) \cap \R^\AA_+$ and hence $-\Omega(\pi)
\cdot \Gamma_\pi$ is an open cone in $H(\pi)$.  In fact, it is the
interior of the cone $H^+(\pi)$ of Veech \cite {V1}.

\begin{lemma} \label {projective}
For every $[h] \in \P H(\pi)$, the set $\Theta(\Pi(\pi)) \cdot
[h]$ contains a dense subset of $\P(-\Omega(\pi) \cdot
\Gamma_\pi)$.
\end{lemma}

\begin{pf}
For almost every $x \in \Upsilon_\RR$, $[h]$ is not symplectically
orthogonal to the span $E^{ss}(x)$ of $\Omega(\pi)\cdot\lambda$.
Consequently, $[h] \not \subset E^c(x) \oplus E^{ss}(x)$. By the
Oseledets theorem, $\widehat B^Z(R_Z^{m-1}(x)) \cdots B^Z(x) \cdot
[h]$ and $\Omega(\pi^{(m)}) \cdot [\tau^{(m)}]=E^{uu}(x^{(m)})$
are asymptotic. On the other hand, by ergodicity, the sequence
$\widehat
R_Z^m(x)=x^{(m)}=(([\lambda^{(m)}],\pi^{(m)}),[\tau^{(m)}])$ is
dense in $\Upsilon_\RR$ for almost every $x$. This implies
$E^{uu}(x^{(m)})$ is dense in $\P(-\Omega(\pi) \cdot \Gamma_\pi)$.
The result follows from these two observations..
\end{pf}

\begin{cor} \label {k=1}
The action, via $\Theta$, of the Rauzy monoid $\Pi(\pi)$ on $\P
H(\pi)$ is minimal.
\end{cor}

\begin{pf}
By Lemma \ref {projective}, the closure of any orbit must intersect every
other orbit.  This implies minimality, obviously.
\end{pf}

\subsection{Choices of permutations in Rauzy classes} \label {perm}

Let us say that $\pi$ is {\it standard} if the first in the
top/bottom is the last in the bottom/top.  Note that a standard
permutation is always irreducible.

\begin{lemma}[Rauzy \cite {R}] \label{Rauzy}
In every Rauzy class there exists a standard permutation.
\end{lemma}


Let $\pi \in \ssigma(\AA)$ be standard. Assuming $d\ge 3$, we call
$\pi$ {\it degenerate} if there exists $B \in \AA$ which is either
second of the top and bottom rows or second to last of the top and
bottom rows. Assuming $d\ge 4$, we call $\pi$ {\it good} if
forgetting the first (and last) letters of the top and the bottom
rows gives an irreducible permutation. Notice that a standard
permutation can not be both degenerate and good.

\begin{lemma}[Lemma 20 of \cite {KZ}]
\label{degenerateorgood} In every Rauzy class with $\#\AA\ge 3$
there exists either a good permutation or a degenerate
permutation.
\end{lemma}

\begin{pf}
We give a proof here for the convenience of the reader, and to
avoid confusion with the slightly different language of \cite
{KZ}. For $d=\# \AA =3$ the result is immediate. In what follows
we suppose $d\geq 4$.

Let $A$ be the first letter in the top and $E$ be the first letter
in the bottom. Suppose, by contradiction, that no standard
permutation in $\RR$ is either degenerate or good. Let $\pi$ be a
standard permutation and $\pi'$ be obtained by forgetting $A$ and
$E$. Then $\pi'$ is reducible, and so there exists a maximal $1
\leq k \leq d-3$ such that the set of first $k$ symbols in the top
and in the bottom of $\pi'$ coincide. Since $\pi$ is
non-degenerate, we must have $2 \leq k \leq d-4$.  Assume $\pi$
has been chosen so that the resulting $k$ is maximal. Let $x_1,
\ldots, x_k$ be the first $k$ letters in the top and $y_1,\ldots,
y_k$ be the first $k$ letters in the bottom of $\pi'$ (so
$\{x_1,\ldots,x_k\}=\{y_1,\ldots,y_k\}$). Let $C$ be the letter in
position $d-1$ in the top of $\pi$, and $l$ be its position in the
bottom of $\pi$. Then $k+2 \leq l \leq d-2$. Let $C'$ be the
letter preceding $C$ in the bottom of $\pi$, and $r$ be its
position in the top of $\pi$. Then $2 \leq r \leq d-2$. Let us
consider the following Rauzy path starting from $\pi$:
\begin{enumerate}
\item Apply $d-r$ bottom iterations to $\pi$,
      so that $C'$ becomes last in the top.
\item Apply $d-l$ top iterations, so that $C$ becomes last in the bottom.
\item Apply $r-1$ bottom iterations, so that $E$ becomes last in the top.
\item Finally, apply $1$ top iteration, so that $A$ becomes last in the bottom.
\end{enumerate}
Notice that step (1) sends $C$ to the position $d-r$ in the top,
preceding $E$, and step (2) sends $A$ to the position $d-1$ in the
bottom. In the end, we get a new standard permutation $\tilde
\pi$: $A$ is last in the bottom and $E$ is last in the top. Let
$\tilde \pi'$ be obtained from $\tilde \pi$ by forgetting the
letters $A$ and $E$. There are two cases to consider.

If $l>k+2$ then $r>k+1$. The calculation for this case is detailed
in the following formula:
{\tiny
\begin{equation*}
\begin{aligned} & \left(\begin{array}{cccccccccccc} A & x_1 &
\cdot & x_k & \cdot & \cdot & \cdot & \cdot & C' & \cdot & C & E
\\
E & y_1 & \cdot & y_k & \cdot & C' & C & \cdot & \cdot & \cdot & \cdot & A
\end{array}\right)
\\ 
& \mapsto
\left(\begin{array}{cccccccccccccccc}
A & \cdot & \cdot & \cdot & \cdot & \cdot & \cdot & \cdot & \cdot & C & E & x_1 & \cdot & x_k & \cdot & C'
\\
E & y_1 & \cdot & y_k & \cdot & C' & C & \cdot & \cdot & \cdot & \cdot & \cdot & \cdot & \cdot & \cdot & A
\end{array}\right)
\\ 
& \mapsto
\left(\begin{array}{cccccccccccccccc}
A & \cdot & \cdot & \cdot & \cdot & \cdot & \cdot & \cdot & \cdot & C & E & x_1 & \cdot & x_k & \cdot & C'
\\
E & y_1 & \cdot & y_k & \cdot & C' & \cdot & \cdot & \cdot & \cdot & \cdot & \cdot & \cdot & \cdot & A & C
\end{array}\right)
\\ 
& \mapsto
\left(\begin{array}{cccccccccccccccc}
A & \cdot & \cdot & \cdot & \cdot & \cdot & \cdot & \cdot & \cdot & C & x_1 & \cdot & x_k & \cdot & C' & E
\\
E & y_1 & \cdot & y_k & \cdot & C' & \cdot & \cdot & \cdot & \cdot & \cdot & \cdot & \cdot & \cdot & A & C
\end{array}\right)
\\ 
& \mapsto
\left(\begin{array}{cccccccccccccccc}
A & \cdot & \cdot & \cdot & \cdot & \cdot & \cdot & \cdot & \cdot & C & x_1 & \cdot & x_k & \cdot & C' & E
\\
E & C & y_1 & \cdot & y_k & \cdot & C' & \cdot & \cdot & \cdot & \cdot & \cdot & \cdot & \cdot & \cdot & A
\end{array}\right)
\end{aligned}
\end{equation*}}
Notice that $C$ precedes $x_1, \ldots, x_k$ in the top and
precedes $y_1, \ldots, y_k$ in the bottom of $\tilde \pi$. By
assumption, $\tilde \pi$ is neither degenerate nor good. The first
of these facts implies that $C$ is not the second letter in the
top of $\tilde\pi$. The second one means that $\tilde \pi'$ is
reducible: there exists $1 \leq \tilde k \leq d-3$ such that the
first $\tilde k$ elements in the top and the bottom of $\tilde
\pi'$ coincide. In view of the previous observations, this implies
that the first $\tilde k$ elements in the bottom of $\tilde \pi'$
include $C$, $y_1$,\ldots,$y_k$. Thus $\tilde k>k$, contradicting
the choice of $k$.

If $l=k+1$, then $C'=y_k=x_{r-1}$. The calculation for this case
is detailed in the formula:
{\tiny
\begin{equation*}
\begin{aligned}
& \left(\begin{array}{cccccccccc}
A & x_1 & \cdot & x_{r-1} & \cdot & x_k & \cdot & \cdot & C & E
\\
E & y_1 & \cdot &\cdot & \cdot & y_k & C & \cdot & D & A
\end{array}\right)
\\ 
& \mapsto
\left(\begin{array}{cccccccccccccc}
A & \cdot & \cdot & D & \cdot & \cdot & \cdot & \cdot & C & E & x_1 & \cdot & \cdot & x_{r-1}
\\
E & y_1 & \cdot &\cdot & \cdot & y_k & C & \cdot & \cdot & \cdot & \cdot & \cdot & D & A\end{array}\right)
\\ 
& \mapsto
\left(\begin{array}{cccccccccccccc}
A & \cdot & \cdot & D & \cdot & \cdot & \cdot & \cdot & C & E & x_1 & \cdot & \cdot & x_{r-1}
\\
E & y_1 & \cdot &\cdot & \cdot & y_k & \cdot & \cdot & \cdot & \cdot & \cdot & D & A & C\end{array}\right)
\\ 
& \mapsto
\left(\begin{array}{cccccccccccccc}
A & \cdot & \cdot & D & \cdot & \cdot & \cdot & \cdot & C & x_1 & \cdot & \cdot & x_{r-1} & E
\\
E & y_1 & \cdot &\cdot & \cdot & y_k & \cdot & \cdot & \cdot & \cdot & \cdot & D & A & C\end{array}\right)
\\ 
& \mapsto
\left(\begin{array}{cccccccccccccc}
A & \cdot & \cdot & D & \cdot & \cdot & \cdot & \cdot & C & x_1 & \cdot & \cdot & x_{r-1} & E
\\
E & C & y_1 & \cdot &\cdot & \cdot & y_k & \cdot & \cdot & \cdot & \cdot & \cdot & D & A \end{array}\right)
\end{aligned}
\end{equation*}}
Let $D$ be the letter in position $d-1$ in the bottom of $\pi$.
Notice that $C \neq D$ since $\pi$ is non-degenerate. After the
first step, $C$ precedes $E$ that precedes
$x_1,\ldots,x_{r-1}=y_k$. In particular, $D$ appears before $C$ in
the top. It follows that $D$ appears before $C$ in the top of
$\tilde \pi$. On the other hand, after the second step, $D$
precedes $A$ in the bottom. It follows that $D$ is in position
$d-1$ in the bottom of $\tilde \pi$. So, $C$ is the first letter
and $D$ is the last letter in the bottom of $\tilde \pi'$, and $D$
appears before $C$ in the top. This implies that $\tilde \pi'$ is
irreducible, which contradicts the hypothesis. So, this case can
not really occur.
\end{pf}

\section{Twisting and pinching} \label {4}

In this section we consider actions of a monoid $\BB$ by linear
isomorphisms of a vector space $H$. We introduce certain
properties of monoids that we call twisting, twisting of isotropic
spaces, pinching, strong pinching, and simplicity, and describe
some logical relations between them.

\subsection{Twisting}

We say that a monoid twists a subspace $F$ if it contains enough
elements to send $F$ outside any finite union of hyperplane sections.
More precisely,

\begin{definition}

We say that $\BB$ {\it twists} some $F \in \grass(k,H)$ if, for every
finite subset $\{F_i\}_{i=1}^m$ of $\grass(\dim H-k,H)$, there exists
$x \in \BB$ such that
\be
(x \cdot F) \cap F_i=\{0\}, \quad 1 \leq i \leq m.
\ee

\end{definition}

In connection with the next lemma, observe that a linear arrangement $S$
is $\BB$-invariant (that is, $x\cdot S=S$ for every $x\in\BB$) if and
only if it is $\BB^{-1}$-invariant, where $\BB^{-1}=\{x^{-1}: x\in\BB\}$.
In particular, the lemma implies that {\it $\BB$ twists $F$ if and only
if $\BB^{-1}$ twists $F$.}

\begin{lemma}\label{twistsiff}

A monoid $\BB$ twists $F \in \grass(k,H)$ if and only if $F$ does not belong
to any non-trivial invariant linear arrangement in $\grass(k,H)$.

\end{lemma}

\begin{pf}

Let $S \subset \grass(k,H)$ be a non-trivial linear arrangement
containing $F$.  Then $S$ is contained in a finite union $\tilde
S$ of hyperplane sections $S_i=\{F' \in \grass(k,H) : F' \cap F_i
\neq \{0\}\}$ for all $1 \leq i \leq m$. If $\BB$ twists $F$ then
there exists $x \in \BB$ such that $(x \cdot F) \cap F_i=\{0\}$,
$1 \leq i \leq m$, that is $x \cdot F \notin \tilde S$. Since
$S\subset\tilde S$, it follows that $S$ is not invariant. On the
other hand, if $\BB$ does not twist $F$ then there exists a finite
union of hyperplane sections $\tilde S \subset \grass(k,H)$ such
that $x \cdot F \in \tilde S$ for every $x \in \BB$. From
Corollary \ref{c.linear9}(1) we get that $S=\cap_{x \in \BB}
x^{-1} \cdot\tilde S$ is a non-trivial linear arrangement
containing $F$. It is clear that $x^{-1} \cdot S \supset S$ for
every $x\in\BB$. From Corollary \ref{c.linear9}(3) it follows that
$x \cdot S=S$ for every $x \in \BB$.
\end{pf}

Next, we prove that given any finite family of hyperplane sections and
subspaces, possibly with variable dimensions, there exists some
isomorphism in $\BB$ that sends every subspace outside the corresponding
hyperplane section (simultaneously), provided $\BB$ twists each one of
the subspaces individually.

\begin{lemma} \label {twist}

For $1\le j \le m$, let $k_j$ satisfy $1 \leq k_j \leq 2g-1$, let
$F_j \in \grass(k_j,H)$, and let $F'_j \in \grass(\dim H-k_j,H)$.
Assume that $F_j$ is twisted by $\BB$ for every $j$.
Then there exists $x \in \BB$ such that $x \cdot F_j \cap F'_j=\{0\}$
for all $1 \leq j \leq m$.

\end{lemma}

\begin{pf}

Let $S_j \subset \Ext^{k_j} H$ be the hyperplane dual to $F'_j$.
Consider the vector space $X=\prod_{j=1}^m \Ext^{k_j} H$ and let
$W=\cap_{x \in \BB} x^{-1} \cdot Y$ where
 \be
Y=\bigcup_{j=1}^m \Big[\prod_{i<j} \Ext^{k_i} H \times S_j \times
\prod_{i>j} \Ext^{k_i} H \Big]\subset X.
 \ee
Let $u_j \in \Ext^{k_j} H$ projectivize to $F_j$. If the
conclusion does not hold then $u=\big(u_j\big)_{j}$ belongs to
$W$. Using Lemmas \ref{l.linear1} and \ref{l.linear3}, we may
write $W$ as a finite union (uniquely defined up to order) of
subspaces $W_l$, $l=1, \ldots, L$ where $L$ is minimal and the
$W_l$ are of the form
 \be
 W_l= W_{l,1} \times \cdots \times W_{l,m},
 \ee
where each $W_{l,j}$ is an intersection of geometric hyperplanes
of $\Ext^{k_j} H$, or else coincides with the whole exterior
product. Given any $x\in\BB$, we have $x^{-1}\cdot W\supset W$ and
so $x^{-1}\cdot W=W$, by Corollary \ref{c.linear9}(3). This means
that each $x \in \BB$ permutes the $W_l$ and so, for every $j$, it
permutes the $W_{l,j}$ (counted with multiplicity). We have $u \in
W_{l_0}$ for some $l_0$, that is, $u_j \in W_{l_0,j}$ for all $j$.
Since $W \neq X$, there exists $j_0$ such that $W_{l_0,j_0} \neq
\Ext^{k_{j_0}} H$, and so $\grass(k_{j_0},H)$ is not contained in
the projectivization of $W_{l_0,j_0}$. Note that, by construction,
the intersection of $\grass(k_{j_0},H)$ with the projectivization
of $W_{l_0,j_0}$ contains $F_{j_0}$. Let $W^0$ be the union of all
$W_{l,j_0}$ whose projectivization intersects but does not contain
$\grass(k_{j_0},H)$. Then the projectivization of $W^0$
intersected with $\grass(k_{j_0},H)$ is a non-trivial linear
arrangement in the Grassmannian, invariant for $\BB$ and
containing $F_{j_0}$. This contradicts the assumption that $\BB$
twists $F_{j_0}$.
\end{pf}

\begin{lemma} \label {minimal implies twisting}

Let $\BB$ be a monoid acting symplectically on $(H,\omega)$, and assume that
the action of $\BB$ on the space of Lagrangian flags $\LL(H)$ is minimal.
Then $\BB$ twists isotropic subspaces.

\end{lemma}

\begin{pf}

If $\BB$ acts minimally on $\LL(H)$ then it also acts minimally on
each $\iso(k,H)$, $1 \leq k \leq g$ (where $\dim H=2g$). We will
only use this latter property. Let $F \in \iso(k,H)$ and $S$ be a
non-trivial invariant linear arrangement in $\grass(k,H)$.  If $F
\in S$ then $S \cap \iso(k,H)$ is a non-empty closed invariant set
under the action of $\BB$. By minimality, it must be the whole of
$\iso(k,H)$.  This contradicts Lemma \ref {2.7}. Therefore,
$F\notin S$. In view of Lemma~\ref{twistsiff}, this proves the
claim.
\end{pf}

\subsection{Pinching}

Assume that $\BB$ acts symplectically on $(H,\omega)$, $\dim H=2g$.

\begin{definition}

We say that $\BB$ is {\it strongly pinching} if for every $C>0$,
there exists $x \in \BB$ such that
\be \label {1111}
\ln \sigma_g(x)>C
\quad\text{and}\quad
\ln \sigma_i(x)>C \ln \sigma_{i+1}(x) \quad\text{for all } 1 \leq i \leq g-1.
\ee
This is independent of the choice of the inner product used to
define the singular values.

\end{definition}

The converse statement in the next lemma says that strong pinching
together with the twisting of all isotropic spaces provide good
separation of the Lyapunov exponents.

\begin{lemma} \label {pinchingsimple}

If for every $C>0$, there exists $x \in \BB$ such that
\be
\theta_g(x)>0
\quad\text{and}\quad
\theta_k(x)>C \theta_{k+1}(x) \quad\text{for all } 1 \leq k \leq g-1,
\ee
then $\BB$ is strongly pinched.  Most important, the converse holds if
$\BB$ twists isotropic spaces.

\end{lemma}

\begin{pf}

Clearly, if $x$ has simple Lyapunov spectrum we have \be
\sigma_k(x^n) \asymp e^{n \theta_k(x)}. \ee This gives the first
assertion. For the second one, let $x_n \in \BB$ be such that \be
\sigma_g(x_n)>n \quad\text{and}\quad \ln \sigma_k(x_n)>n \ln
\sigma_{k+1}(x_n) \quad\text{for all } 1 \leq k \leq g-1. \ee We
may assume that $x_n \cdot E^+_k(x_n)$ converges to some $E^u_k
\in \grass(k,H)$ and $E^-_k(B_n)$ converges to some $E^s_k \in
\grass(2g-k,H)$, $1 \leq k \leq g$. By Lemma~\ref{isolimit}, the
subspaces $E^u_k$ are isotropic. It follows from Lemma \ref{twist}
that if $\BB$ twists isotropic subspaces, there exists $x \in \BB$
such that \be x \cdot E^u_k \cap E^s_k=\{0\} \quad\text{for all }
1 \leq k \leq g. \ee By Lemma~\ref {sigma theta}, there exists
$C>0$ such that \be |\theta_k(x x_n)-\ln \sigma_k(x_n)|<C
\quad\text{for all } 1 \leq k \leq g, \ee which implies the
result.
\end{pf}

The next lemma indicates a useful special situation where one has the
strong pinching property.

\begin{lemma} \label {criterion}

If for every $C>0$ there exists $x \in \BB$ for which $1$ is an
eigenvalue of geometric multiplicity $1$ (dimension of the
eigenspace equal to $1$), and we have \be \theta_{g-1}(x)>0
\quad\text{and}\quad \theta_k(x)>C\theta_{k+1}(x) \quad\text{for
all } 1 \leq k \leq g-2, \ee then $\BB$ is strongly pinched.

\end{lemma}

\begin{pf}

Since the action is symplectic, the eigenvalue $1$ has even
algebraic multiplicity. Considering the Jordan form of $x$, the
hypotheses imply
 \be
 \sigma_g(x^n) \asymp n
 \ee
 \be
 \sigma_k(x^n) \asymp e^{\theta_k(x) n} \quad\text{for all } 1
 \leq k \leq g-1.
 \ee
The claim follows immediately.
\end{pf}

\begin{lemma} \label {more}

Assume $\BB_0 \subset \BB$ twists isotropic subspaces and is
strongly pinching. Then for every $x \in \BB$ and any $C>0$ there
exists $x_0 \in \BB_0$ such that \be \label{primeira}
\theta_g(x_0)>0 \quad\text{and}\quad \theta_i(x_0)>C
\theta_{i+1}(x_0) \quad\text{for all } 1 \leq i \leq g-1, \ee \be
\label{segunda} \theta_g(x x_0)>0 \quad\text{and}\quad \theta_i(x
x_0)>C\theta_{i+1}(x x_0) \quad\text{for all } 1 \leq i \leq g-1.
\ee

\end{lemma}

\begin{pf}

Let $x_2 \in \BB_0$ satisfy the conditions in \eqref{primeira}:
$\theta_g(x_2)>0$ and $\theta_i(x_2)>C \theta_{i+1}(x_2)$ for all
$1 \leq i \leq g-1$. Let $E^u_k(x_2)=\lim x^n_2 \cdot
E^+_k(x^n_2)$, which is the space spanned by the $k$ eigenvectors
with largest eigenvalues, and let $E^s_k(x_2)=\lim E^-_k(x^n_2)$
which is the space spanned by the $2g-k$ remaining eigenvectors.
Using Lemma \ref{twist}, select $x_1 \in \BB_0$ such that $x_1
\cdot E^u_k \cap E^s_k=\{0\}$ and $x x_1 \cdot E^u_k \cap
E^s_k=\{0\}$ for all $1 \leq k \leq g$.  Then $x_0=x_1 x^n_2$
satisfies all the conditions for $n$ large, by the same argument
as in Lemma~\ref{pinchingsimple}.
\end{pf}

\subsection{Simple actions}

\begin{definition}

We say that $\BB$ is {\it twisting} if it twists any $F \in
\grass(k,H)$ for any $1 \leq k \leq \dim H-1$.

\end{definition}

From the observation preceding Lemma~\ref{twistsiff} we have that
$\BB$ is twisting if and only if $\BB^{-1}$ is twisting.

\begin{definition}

We say that $\BB$ is {\it pinching} if for every $C>0$ there
exists $x \in \BB$ such that \be \sigma_i(x)>C \sigma_{i+1}(x)
\quad\text{for all } 1 \leq i \leq \dim H-1. \ee This is
independent of the choice of the inner product used to define the
singular values.

\end{definition}

\begin{definition}

We say that $\BB$ is {\it simple} if it is both twisting and pinching.

\end{definition}

\begin{lemma} \label {b-1}
Let $\BB$ be a simple monoid.  Then the inverse monoid $\BB^{-1}$
is also simple.
\end{lemma}

\begin{pf}
Twisting follows directly from the observation preceding Lemma
\ref{twistsiff}. Pinching follows from the fact that
$\ln\sigma_i(x)=-\ln\sigma_{\dim H-i+1}(x^{-1})$ for all $1 \leq i
\leq \dim H-1$.
\end{pf}

We suspect there could be sufficient conditions for twisting along
the following lines:

\begin{problem}

If $\BB$ acts minimally on $\P H$ is it necessarily twisting? For
symplectic actions, one can ask a weaker question: Does minimal
action on $\P H$ imply twisting of isotropic subspaces?

\end{problem}

\begin{lemma} \label {tps}

Let $\BB$ be a monoid acting symplectically on $(H,\omega)$.  If $\BB$
twists isotropic subspaces and is strongly pinching then it is simple.

\end{lemma}

\begin{pf}

To see that strong pinching implies pinching, it is enough to
consider an adapted inner product, for which \eqref {1111} implies
also \be -\ln \sigma_{g+1}(x)>C, \ee \be -\ln \sigma_{i+1}(x)>-C
\ln \sigma_i(x) \quad\text{for all}\quad g+1 \leq i \leq 2g-1, \ee
where $\dim H=2g$. We are left to show that, under the hypotheses,
$\BB$ twists any $F \in \grass(k,H)$.  It is easy to see that
$\BB$ twists $F$ if and only if it twists its symplectic
orthogonal. So it is enough to consider the case $1 \leq k \leq
g$.  Let $S \subset \grass(k,H)$ be a non-trivial invariant linear
arrangement.  By Lemma \ref {pinchingsimple}, there exists $x_1
\in \BB$ with simple Lyapunov spectrum.  Let $E^s=\lim
E^-_g(x_1^n)$ and $E^u=\lim x_1^n \cdot E^+_g(x_1^n)$ be the
stable and unstable eigenspaces of $x_1$.  Thus both $E^s$ and
$E^u$ are isotropic subspaces of $H$.  If $F \in S$ then there
exists $x_0 \in \BB$ such that $(x_0 \cdot F) \cap E^s=\{0\}$.  It
follows that $x_1^n x_0 \cdot F$ converges to a subspace of $E^u$.
Since $S$ is closed, we conclude that $S \cap \iso(k,H) \neq
\emptyset$. This contradicts the assumption that $\BB$ twists all
elements of $\iso(k,H)$.
\end{pf}

\begin{rem} \label {similar}

An argument similar to the proof of the previous lemma shows that if $\BB$
is simple and acts symplectically then any closed invariant set in the flag
space $\FF(H)$ intersects the embedding of the space of
Lagrangian flags $\LL(H)$, that is, it contains some $(F_i)_{i=1}^{2g-1}$
such that $(F_i)_{i=1}^g$ is a Lagrangian flag and $F_{2g-i}$ is the
symplectic orthogonal of $F_i$.

\end{rem}

\begin{lemma} \label {large}

Let $\BB_0 \subset \BB$ be a large submonoid in the sense that
there exists a finite subset $Y \subset \BB$ and $z \in \BB$
such that for every $x \in \BB$, $y x z \in \BB_0$ for some
$y \in Y$.  If $\BB$ is twisting or pinching then $\BB_0$ also is.
Assuming the action of $\BB$ is symplectic, if $\BB$ twists isotropic
subspaces or is strongly pinching then the same holds for $\BB_0$.

\end{lemma}

\begin{pf}

Notice that $|\ln \sigma_i(x)-\ln \sigma_i(yxz)|<C$ where $C$ only
depends on $y$, $z$ and the choice of the inner product used to
define the singular values.  Thus if $\BB$ is (strongly) pinching
then $\BB_0$ is also (strongly) pinching. Let $S \subset
\grass(k,H)$ be a non-trivial linear arrangement invariant for
$\BB_0$. Then
 \be
 \bigcap_{x \in \BB} x^{-1} \cdot \bigcup_{y \in Y} y^{-1} \cdot S
 \supset \bigcap_{y \in Y} \bigcap_{x_0 \in \BB_0} (y^{-1} x_0 z^{-1})^{-1}
         \cdot \bigcup_{y \in Y} y^{-1} \cdot S
 \supset z \cdot \bigcap_{y \in Y} \bigcap_{x_0 \in \BB_0} x_0^{-1} \cdot S = z \cdot S.
\ee
So, $z \cdot S$ is contained in a non-trivial linear arrangement invariant
for $\BB$. This shows that if $\BB$ is twisting then so is $\BB_0$.
If $\BB$ acts symplectically and $S$ intersects $\iso(k,H)$ then $z \cdot S$
also does, so if $\BB$ twists isotropic subspaces then $\BB_0$ also does.
\end{pf}

\begin{lemma}
Let $\BB$ be a simple monoid.  Then there exists $x \in \BB$ with simple
Lyapunov spectrum $\theta_i(x)>\theta_{i+1}(x)$, $1 \leq i \leq \dim H-1$.
\end{lemma}

\begin{pf}
Let $x_n \in \BB$ be such that $\ln \sigma_i(x)-\ln \sigma_{i+1}(x)>n$,
$1 \leq i \leq \dim H-1$.  We may assume that $E^-_i(x_n)$ and $x_n \cdot
E^+_i(x_n)$ converge to spaces $E^s_i$ and $E^u_i$, $1 \leq i \leq \dim H-1$.
Let $x \in \BB$ be such that $(x \cdot E^u_i) \cap E^s_i=\{0\}$,
$1 \leq i \leq \dim H-1$.  By Lemma~\ref {sigma theta} there exists
$C>0$ such that $|\theta_i(x x_n)-\ln \sigma_i(x_n)|<C$,
$1 \leq i \leq \dim H-1$.  This gives the result.
\end{pf}

\section{Twisting of Rauzy monoids} \label {5}

Recall that a Rauzy monoid $\Pi(\pi)$ acts symplectically on $H(\pi)$
(by $\gamma \cdot h=\Theta(\gamma) \cdot h$).
In particular, it acts on the space of Lagrangian flags $\LL(H(\pi))$.
Our aim in this section is to prove the following result:

\begin{thm} \label {minimal action}
Let $\pi$ be irreducible. The action of the Rauzy monoid
$\Pi(\pi)$ on the space of Lagrangian flags $\LL(H(\pi))$ is
minimal.
\end{thm}

\begin{cor} \label {twisting}
Let $\pi$ be irreducible. The action of the Rauzy monoid
$\Pi(\pi)$ on $H(\pi)$ twists isotropic subspaces.
\end{cor}

\begin{pf}
This follows directly from Theorem \ref {minimal action} and
Lemma \ref {minimal implies twisting}.
\end{pf}

%
%
%
%
%
%

\subsection{Simple reduction}

Let $\AA$ be an alphabet on $d \geq 3$ symbols, $B \in \AA$, and
$\AA'=\AA \setminus \{B\}$.  Given $\pi \in \ssigma(\AA)$, let
$\pi'$ be obtained from $\pi$ by removing $B$ from the top and
bottom rows.  If $\pi' \in \ssigma(\AA')$, we say that $\pi'$ is a
{\it simple reduction} of $\pi$. Let $2g(\pi)=\dim
H(\pi)=\rank\Omega(\pi)$ and analogously for $\pi'$. Let $P:\R^\AA
\to \R^{\AA'}$ be the natural projection, and $P^*:\R^{\AA'}\to
\R^{\AA}$ be its adjoint (the natural inclusion). Notice that \be
\label {pompstar} P \Omega(\pi) P^*=\Omega(\pi'). \ee

\begin{lemma} \label {simplereduction}

Let $\pi'$ be a simple reduction of $\pi$.  Then either $g(\pi)=g(\pi')$ or
$g(\pi)=g(\pi')+1$.  Moreover, the following are equivalent:
\begin{enumerate}
\item $g(\pi)=g(\pi')$,
\item $H(\pi)$ is spanned by $\{\Omega(\pi) \cdot e_x,\ x \in \AA'\}$,
\item $e_B \notin H(\pi)$,
\item $e_B$ does not belong to the span of
$\{\Omega(\pi) \cdot e_x,\ x \in \AA'\}$,
\item $P$ restricts to a symplectic isomorphism $H(\pi) \to H(\pi')$.
\end{enumerate}

\end{lemma}

\begin{pf}

There are two possibilities for the value of $\rank P \Omega(\pi) P^*$:
\begin{enumerate}
\item [(a)] We may have
$\rank P \Omega(\pi) P^*=\rank \Omega(\pi)$. Since the rank can not
increase by composition, this implies that
$\rank P\Omega(\pi) P^*=\rank P \Omega(\pi)=\rank \Omega(\pi) P^*=\rank
\Omega(\pi)$.
\item [(b)] We may have $\rank P\Omega(\pi)P^*<\rank \Omega(\pi)$.
Notice that $\rank P\Omega(\pi)P^* \geq \rank \Omega(\pi)-2$.
Since $\rank \Omega(\pi)$ and $\rank P\Omega(\pi)P^*$ are even, we must
have $\rank P\Omega(\pi)P^*=\rank \Omega(\pi)-2$.  This implies that $\rank
P\Omega(\pi)=\rank \Omega(\pi)P^*=\rank \Omega(\pi)-1$.
\end{enumerate}
By \eqref {pompstar}, in case (a) we have $g(\pi)=g(\pi')$ and in case (b)
we have $g(\pi)=g(\pi')+1$. Notice that
\begin{enumerate}
\item is equivalent to
      $\rank P\Omega(\pi) P^*=\rank \Omega(\pi)$,
\item is equivalent to
      $\Omega(\pi) \cdot \R^\AA=\Omega(\pi) (P^* \cdot \R^{\AA'})$,
      that is, $\rank\Omega(\pi)=\rank \Omega(\pi) P^*$,
\item is equivalent to
      $(\Omega(\pi) \cdot \R^\AA) \cap \Ker P=\{0\}$,
      that is, $\rank P \Omega(\pi)=\rank \Omega(\pi)$,
\item is equivalent to
      $(\Omega(\pi) (P^* \cdot \R^{\AA'}))\cap \Ker P=\{0\}$,
      that is $\rank P \Omega(\pi) P^*=\rank \Omega(\pi) P^*$.
\end{enumerate}
Thus, (1), (2), (3), (4) are all equivalent.
It is clear that (5) implies (1).
Notice that, by \eqref{pompstar}, we have $P \cdot H(\pi) \supset H(\pi')$.
So, (1) implies that $P$ restricts to an isomorphism $H(\pi) \to H(\pi')$.
To see that this isomorphism is symplectic, let $\omega$ and $\omega'$ be the
symplectic forms on $H(\pi)$ and $H(\pi')$.  Then, using \eqref{pompstar}
and the definition of $\omega$ and $\omega'$,
$$
\begin{aligned}
\omega'(P \Omega(\pi) \cdot e_x,P \Omega(\pi) \cdot e_y)
& =\omega'(\Omega(\pi') \cdot e_x,P\Omega(\pi) \cdot e_y \rangle
=\langle e_x, P\Omega(\pi) \cdot e_y \rangle
\\ & =\langle P^*\cdot e_x, \Omega(\pi) \cdot e_y \rangle
=\omega(\Omega(\pi) \cdot e_x,\Omega(\pi) \cdot e_y)
\end{aligned}
$$
for all $x, y \in \AA'$. In view of (2), this implies that
$P:H(\pi) \to H(\pi')$ is symplectic.
\end{pf}

\subsection{Simple extension}

Let $\AA'$ be an alphabet on $d \geq 2$ symbols and $\pi' \in
\ssigma(\AA')$. Let $A$ be the first in the top and $E$ be the
first in the bottom. If $B \notin \AA'$, let $\AA=\AA' \cup
\{B\}$. Let $C,D \in \AA'$ be such that $(A,E) \neq (C,D)$ (we
allow $C=D$). Let $\pi=\LE(\pi')$ be obtained by inserting $B$ in
$\pi'$ just before $C$ in the top and before $D$ in the bottom.
This operation is described by:
{\tiny
$$
\pi'=\left(\begin{array}{cccccc}
     A & \cdot & C & \cdot & \cdot & \cdot \\
     E & \cdot & \cdot & \cdot & D & \cdot
\end{array} \right) \mapsto
\pi =\left(\begin{array}{cccccccc}
     A & \cdot & B & C & \cdot & \cdot & \cdot & \cdot \\
     E & \cdot & \cdot & \cdot & \cdot & B & D & \cdot
\end{array}\right)
$$}
\begin{lemma}
$\pi=\LE(\pi')$ is irreducible.
\end{lemma}

\begin{pf}
Suppose the first $k<\#\AA$ symbols in the top and the bottom of
$\pi$ coincide. Note that $k\ge 2$. If these symbols do not
include $B$, then they are also the first $k$ symbols in the top
and the bottom of $\pi'$. Otherwise, removing $B$ we obtain the
first $k-1$ symbols in the top and the bottom of $\pi'$. In either
case, this contradicts the assumption that $\pi'$ is irreducible.
\end{pf}

This immediately extends to a map $\LE$ defined (by the same rule)
on the whole Rauzy class $\RR(\pi')$. We are going to see that
$\LE$ takes values in the Rauzy class of $\pi$. Given an arrow
$\gamma'$ in $\RR(\pi')$, let $\gamma=\LE_*(\gamma')$ be the path
defined as follows:
\begin{enumerate}
\item If $C$ is last in the top of $\pi'$ and $\gamma'$ is a
bottom arrow then $\gamma$ is the sequence of the two bottom
arrows starting at the $\LE$-image of the start of $\gamma'$. This
is described by:
{\tiny
$$
\begin{aligned}
\gamma':
  \left(\begin{array}{cccccc} \cdot & * & \cdot & \cdot & \cdot & C
\\
\cdot & \cdot & \cdot & D & \cdot & *
\end{array}\right)
& \mapsto
\left(\begin{array}{ccccccc}
\cdot & * & C & \cdot & \cdot & \cdot & \cdot
\\
\cdot & \cdot & \cdot & \cdot & D & \cdot & *
\end{array}\right)
\\
\LE_*(\gamma'):
\left(\begin{array}{cccccccc}
\cdot & * & \cdot & \cdot & \cdot & \cdot & B & C
\\
\cdot & \cdot & \cdot & B & D & \cdot & \cdot & *
\end{array}\right)
& \mapsto
\left(\begin{array}{ccccccccc}
\cdot & * & C & \cdot & \cdot & \cdot & \cdot & B
\\
\cdot & \cdot & \cdot & \cdot & B & D & \cdot & *
\end{array}\right)
\\ & \mapsto
\left(\begin{array}{cccccccccc}
\cdot & * & B & C & \cdot & \cdot & \cdot & \cdot & \cdot
\\
\cdot & \cdot & \cdot & \cdot & \cdot & B & D & \cdot & *
\end{array}\right)
\end{aligned}
$$}

\item If $D$ is last in the bottom of $\pi'$ and $\gamma'$ is a
top arrow then $\gamma$ is the sequence of two bottom arrows
starting at the $\LE$-image of the start of $\gamma'$. This is
analogous to the previous case.

\item Otherwise, $\gamma$ is the arrow of the same type as
$\gamma'$ starting at the $\LE$-image of the start of $\gamma'$.
For example:
{\tiny
$$
\begin{aligned}
\gamma':
\left(\begin{array}{ccccccccc}
\cdot & * & \cdot & \cdot & \cdot & C & \cdot & \cdot & \cdot
\\
\cdot & \cdot & \cdot & D & \cdot & \cdot & \cdot & \cdot & *
\end{array}\right)
& \mapsto
\left(\begin{array}{ccccccccc}
\cdot & * & \cdot & \cdot & \cdot & \cdot & C & \cdot & \cdot
\\
\cdot & \cdot & \cdot & D & \cdot & \cdot & \cdot & \cdot & *
\end{array}\right)
\\
\LE_*(\gamma'):
\left(\begin{array}{ccccccccccc}
\cdot & * & \cdot & \cdot & \cdot & \cdot & B & C & \cdot & \cdot & \cdot
\\
\cdot & \cdot & \cdot & B & D & \cdot & \cdot & \cdot & \cdot & \cdot & *
\end{array}\right)
& \mapsto
\left(\begin{array}{ccccccccccc}
\cdot & * & \cdot & \cdot & \cdot & \cdot & \cdot & B & C & \cdot & \cdot
\\
\cdot & \cdot & \cdot & B & D & \cdot & \cdot & \cdot & \cdot & \cdot & *
\end{array}\right)
\end{aligned}
$$}
\end{enumerate}
In all cases, $\gamma$ starts at the image by $\LE$ of the start of $\gamma'$
and ends at the image by $\LE$ of the end of $\gamma'$.
This shows that $\gamma$ is a path in the Rauzy class of $\pi$, and that
$\LE$ takes values in $\RR(\pi)$.
We extend $\LE_*$ to a map $\Pi(\RR(\pi')) \to \Pi(\RR(\pi))$ in the whole
space of paths, compatible with concatenation, and call
$\LE:\RR(\pi') \to \RR(\pi)$ and $\LE_*:\Pi(\RR(\pi')) \to \Pi(\RR(\pi))$
{\it extension maps}.

\begin{rem}\label{converse}

If $\pi$ is a simple extension of $\pi'$ then $\pi'$ is a simple reduction
of $\pi$. The converse is true if and only if the omitted letter is not
the last on the top nor on the bottom of $\pi$.

\end{rem}

Observe that if $\gamma'=\EE_*(\gamma)$ then
 \be \label {Ptheta}
 P\Theta(\gamma)=\Theta(\gamma') P.
 \ee
Indeed, this may be rewritten as follows (recall that $P(e_B)=0$
and $P(e_x)=e_x$ for every $x\neq B$):
 \be \label{omegamega}
 \langle \Theta(\gamma) \cdot e_x,e_y \rangle
 = \langle \Theta(\gamma') \cdot P(e_x),e_y \rangle
 \quad\text{for all $x \in \AA$ and $y \in \AA'$.}
 \ee
It is enough to check the case when $\gamma'$ is an arrow, because
$\LE_*$ is compatible with concatenation. In case (1) of the
definition above, $\Theta(\gamma')\cdot e_* = e_C +e_*$ and
$\Theta(\gamma')\cdot e_x=e_x$ for any $x\in\AA'\setminus\{*\}$.
On the other hand, $\Theta(\gamma)\cdot e_*=e_B+e_C+e_*$ and
$\Theta(\gamma')\cdot e_x=e_x$ for any $x\in\AA\setminus\{*\}$. In
particular, $\Theta(\gamma)\cdot e_B = e_B$. The claim
\eqref{Ptheta} follows in this case, and the other two are
analogous.

\begin{lemma} \label {g=g'}

Let $\pi$ be a simple extension of $\pi'$, and assume that $g(\pi)=g(\pi')$.
There exists a symplectic isomorphism $H(\pi) \to H(\pi')$ which conjugates
the actions of $\LE_*(\gamma')$ on $H(\pi)$ and of
$\gamma'$ on $H(\pi')$, for every $\gamma' \in \Pi(\pi')$.

\end{lemma}

\begin{pf}

By Lemma \ref{simplereduction}, the natural projection $P:\R^\AA
\to \R^{\AA'}$ restricts to a symplectic isomorphism $H(\pi) \to
H(\pi')$. Then \eqref{Ptheta} shows this isomorphism conjugates
$\Theta(\gamma)\mid H(\pi)$ and $\Theta(\gamma')\mid H(\pi')$.
\end{pf}

\begin{lemma} \label {sex}

Let $\pi \in \ssigma(\AA)$.  If $\# \AA \geq 3$ then there exists $B \in
\AA$ and $\pi' \in \ssigma(\AA \setminus \{B\})$ such that $\pi$ is a simple
extension of $\pi'$.

\end{lemma}

\begin{pf}

Let $A$ be the first in the top and $E$ be the first in the bottom
of $\pi$. If $\pi_0(E)<\pi_1(A)$, let $B=E$. If
$\pi_1(A)<\pi_0(E)$, let $B=A$. If $\pi_0(E)=\pi_1(A)<\# \AA$, let
$B \in \{A,E\}$ be arbitrary. If $\pi_0(E)=\pi_1(A)=\# \AA$, let
$B \in \AA \setminus \{A,E\}$ be arbitrary. Let $\pi'$ be obtained
by forgetting $B$ on $\pi$. Notice that $\pi'$ is irreducible.
Indeed, suppose the first $k$ letters in the top coincide with the
first $k$ letters in the bottom, $1 \le k \le d-2$. In the first
case, one must have $k\ge\pi_1'(A)$ and then, adding the letter
$E$ to the list, one gets that, for $\pi$, the first $k+1$ letters
in the top coincide with the first $k+1$ letters in the bottom.
This contradicts the assumption that $\pi$ is irreducible. A
symmetric argument applies in the second case, and the same
reasoning applies also in the third case. Finally, admissibility
is obvious in the fourth case, because the permutation is
standard. In all cases, $\pi$ is a simple extension of $\pi'$.
\end{pf}

\subsection{Proof of Theorem~\ref {minimal action}}

We prove the theorem by induction on the size $\#\AA$ of the alphabet.
By Corollary \ref {k=1}, the conclusion holds when $\#\AA=2$, since
in this case $\LL(H(\pi))$ is just $\P H(\pi)$.
Assume that it holds for $\#\AA=d-1 \geq 2$, and let us show that it
also does for $\# \AA=d$.

By Lemma~\ref {sex}, there exists $\pi' \in \ssigma(\AA')$,
$\AA'=\AA \setminus \{B\}$ such that $\pi$ is a simple extension
of $\pi'$. If $g(\pi)=g(\pi')$, the result follows immediately
from the induction hypothesis and Lemma~\ref {g=g'}. Thus, we may
assume $g(\pi)=g(\pi')+1$. Let $H=H(\pi)$, $H'=H(\pi')$,
$\Omega=\Omega(\pi)$, $\Omega'=\Omega'(\pi)$, $\omega=\omega_\pi$,
$\omega'=\omega_{\pi'}$. If $\gamma' \in \Pi(\pi')$ then
$\gamma=\LE_*(\gamma')$ is contained in the stabilizer of $e_B$,
that is
 \be
 \Theta(\gamma) \cdot e_B=e_B
 \ee
(because $B$ never wins).

\begin{lemma} \label {k>1}

There is a symplectic isomorphism $H^{\lambda_B} \to H'$ that conjugates
the action of $\LE_*(\gamma')$ on $H^{\lambda_B}$ and the action of
$\gamma'$ on $H'$, for every $\gamma' \in \Pi(\pi')$.

\end{lemma}

\begin{pf}

Write $h \in H$ as $h=\sum_{x \in \AA} u_x (\Omega \cdot e_x)$
with $u_x\in\R$. Then the condition $\omega(h,e_B)=0$ corresponds
to
 \be
 \omega(h,e_B)=\sum_{x \in \AA} u_x \langle e_x,e_B \rangle
 = u_B =0.
 \ee
In other words, $h\in H_{\lambda_B}$ if and only $h$ belongs to
the span of $\Omega \cdot e_x$, $x \in \AA'$. Using
\eqref{pompstar} we obtain that $P:\R^\AA\to\R^{\AA'}$ takes
$H_{\lambda_B}$ to $H'$. The quotient $P:H^{\lambda_B} \to H'$ is
symplectic: using again \eqref{pompstar}, the definition of
$\omega$ and $\omega'$, and \eqref{omegamega},
$$
\omega'(P\Omega \cdot e_x,P\Omega \cdot e_y)
=\omega'(\Omega' \cdot e_x,\Omega' \cdot e_y)
=\langle e_x,\Omega' \cdot e_y \rangle
=\langle e_x,\Omega \cdot e_y \rangle
= \omega(\Omega \cdot e_x,\Omega \cdot e_y)
$$
for all $x, y \in \AA'$.
Moreover, by \eqref{Ptheta}, this map conjugates the action of
$\gamma$ on $H^{\lambda_B}$ with the action of $\gamma'$ on $H'$.
\end{pf}

By the induction hypothesis, we conclude that $\LE_*(\Pi(\pi')) \subset
\Pi(\pi)$ acts on minimally $\LL(H^{\lambda_B})$.
By Lemma~\ref{2.2} and Corollary~\ref{k=1}, this implies that the action
of $\Pi(\pi)$ on $\LL(H)$ is minimal.
\qed.

%
%
%
%

\section{Pinching of Rauzy monoids} \label {pinra}

Our aim in this section is to prove the following result.

\begin{thm} \label {strongpinching}

Let $\pi$ be irreducible. The action of the Rauzy monoid
$\Pi(\pi)$ on $H(\pi)$ is strongly pinching.

\end{thm}

\begin{cor} \label {simple}

Let $\pi$ be irreducible. The action of the Rauzy monoid
$\Pi(\pi)$ on $H(\pi)$ is simple.

\end{cor}

\begin{pf}

This follows directly from Corollary \ref {twisting},
Theorem \ref {strongpinching}, and Lemma \ref {tps}.
\end{pf}

\begin{rem}

In view of Corollary \ref {simple}, we may apply Remark \ref
{similar} to the action of the Rauzy monoid on the space of flags:
any closed invariant set in $\FF(H(\pi))$ intersects the embedding
of the space $\LL(H(\pi))$ of Lagrangian flags. Since the action
on the latter is minimal, by Theorem \ref{minimal action}, we get
that the only minimal set of the action of $\Pi(\pi)$ on the space
of flags is the embedding of Lagrangian flags.

\end{rem}

We start the proof of Theorem~\ref{strongpinching} by observing that
the strong pinching property only depends on the Rauzy class:

\begin{lemma} \label {pinchingrauzy}

If $\pi$ is such that the action of $\Pi(\pi)$ on $H(\pi)$ is strongly
pinching and $\tilde \pi \in \RR(\pi)$ then the action of
$\Pi(\tilde \pi)$ on $H(\tilde \pi)$ is strongly pinching.

\end{lemma}

\begin{pf}

Let $\gamma_0 \in \Pi(\RR(\pi))$ be a path starting at $\tilde
\pi$ and ending at $\pi$, and let $\gamma_1 \in \Pi(\RR(\pi))$ be
a path starting at $\pi$ and ending at $\tilde \pi$.  Then
$\Pi(\tilde \pi) \supset \gamma_0 \Pi(\pi) \gamma_1$, and we
conclude as in the first part of the proof of Lemma~\ref {large}
that $\Pi(\tilde \pi)$ is strongly pinching.
%
\end{pf}

\subsection{Minimal Rauzy classes}

Let us call a Rauzy class $\RR \subset \ssigma(\AA)$ minimal if
$\# \AA=2g(\RR)$.  Recall the definitions of degenerate
permutation and good permutation in Section \ref {perm}. In
particular, these permutations are standard, by definition.

\begin{lemma}\label{degenerate}

Let $\pi \in \ssigma(\AA)$ be a degenerate permutation.
Then $e_B \notin H(\pi)$, where $B$ is the letter appearing in
the second (or second to last) position of both top and bottom.

\end{lemma}

\begin{pf}

Let $\Omega=\Omega(\pi)$. Let $A$ be first in the top (last in the
bottom) and $E$ be the first in the bottom (last in the top).
If $e_B \in H(\pi)$ then, by the equivalence of (3) and (4) in
Lemma \ref {simplereduction}, we can write
 \be
e_B=\sum_{x \neq B} u_x(\Omega \cdot e_x) \quad\text{with }
u_x\in\R.
 \ee
Since $\pi$ is standard, with first/last letters $A$ and $E$, the
definition \eqref{defOmega} gives $\langle \Omega\cdot e_x, e_A
\rangle = 1$ for all $x\neq A$, and $\langle \Omega\cdot e_x, e_E
\rangle = -1$ for all $x\neq E$. Thus, the previous relation
implies
 \be
 0=\langle e_B,e_A \rangle=\sum_{x \neq A,B} u_x
 \quad\text{and}\quad
 0=\langle e_B,e_E \rangle=-\sum_{x \neq E,B} u_x.
 \ee
This implies that $u_A-u_E=0$. On the other hand, \be 1 = \langle
e_B,e_B \rangle = \sum_{x\neq B} \langle \Omega, e_B \rangle = u_E
- u_A \ee because $B$ is the second letter in both top and bottom
(use \eqref{defOmega} once more). This contradicts the previous
conclusion, and this contradiction proves that $e_B\notin H(\pi)$.
\end{pf}

\begin{lemma} \label {goodexists}

Any minimal Rauzy class $\RR \subset \ssigma(\AA)$
with $g(\RR) \geq 2$ contains a good permutation.

\end{lemma}

\begin{pf}

By Lemma \ref {degenerateorgood}, if $\RR$ does not contain a good
permutation then it contains some degenerate permutation $\pi$.
Let $B \in \AA$ be as in Lemma~\ref{degenerate}, and let $\pi'$ be
obtained by forgetting $B$ in $\pi$. Then $\pi'$ is a simple
reduction of $\pi$, and by the equivalence of (1) and (3) in Lemma
\ref {simplereduction}, $g(\pi')=g(\pi)$. Thus $\# \AA \geq
2g(\pi')+1=2g(\pi)+1$, and so $\RR$ is not minimal.
\end{pf}

\begin{lemma} \label {211}

If $\pi$ is good and $\pi''$ is obtained by forgetting the first (and last)
letters of the top and bottom rows then $g(\pi) \leq g(\pi'')+1$.
In particular, if $\RR(\pi)$ is minimal then $\RR(\pi'')$ is also minimal.

\end{lemma}

\begin{pf}

To prove the first claim, let $A$ be first in the top (last in the
bottom) and $E$ be first in the bottom (last in the top) of $\pi$.
Let $\pi'$ be obtained by forgetting $A$ in $\pi$. Then $\pi''$ is
a simple reduction of $\pi'$, which is itself a simple reduction
of $\pi$. Let $\Omega=\Omega(\pi)$, $\Omega'=\Omega(\pi')$,
$\Omega''=\Omega(\pi'')$, $g=g(\pi)$, $g'=g(\pi')$, and
$g''=g(\pi'')$. If $g=g'$ then $g=g' \leq g''+1$, by Lemma~\ref
{simplereduction}, and the conclusion follows. Otherwise, $g=g'+1$
and, by the equivalence of (1) and (4) in Lemma \ref
{simplereduction}, we can write
 \be
 e_A=\sum_{x \neq A} u_x (\Omega \cdot e_x)
 \quad\text{with } u_x\in\R.
 \ee
This implies that (recall the definition \eqref{defOmega} of $\Omega$)
 \be
0 = \langle e_A,e_E \rangle = -\sum_{x \neq A,E} u_x,
\quad\text{and}\quad
1 = \langle e_A,e_A \rangle = \sum_{x \neq A} u_x
 \ee
so that $u_E=1$. From \eqref{defOmega} we have that
$\Omega'\cdot e_x=\Omega\cdot e_x - e_A$ for all $x\neq A$.
Consequently, using the previous equalities,
$$
\Omega'\cdot e_E
= \Omega\cdot e_E - \sum_{x\neq A} u_x (\Omega \cdot e_x) = -
\sum_{x\neq A,E} u_x(\Omega \cdot e_x)
= - \sum_{x\neq A, E} u_x(\Omega' \cdot e_x)
$$
By the equivalence of (1) to (2) in Lemma~\ref{simplereduction},
this gives that $g'=g''$ and so $g=g''+1$.

The last claim in the lemma is an immediate consequence of the
first one, because $\#\AA=2 g(\pi)$ and
$\#\AA-2=\#(\AA\setminus\{A,E\})\ge 2g(\pi'')$, and so the
equality must hold if $g(\pi)\le g(\pi'')+1$.
\end{pf}

\subsection{Reduction to the case of minimal Rauzy classes}

\begin{lemma} \label {dg}

Let $\RR \subset \ssigma(\AA)$ be a Rauzy class such that $\# \AA>2g(\RR)$.
Then there
exists $\pi \in \RR$ which is a simple extension of some $\pi'$ with
$g(\pi)=g(\pi')$.

\end{lemma}

\begin{pf}

By Lemma \ref{Rauzy}, we may consider some standard permutation
$\tilde \pi \in \RR$. Since $\# \AA>2g$, there exists $B \in \AA$
such that $e_B \notin H(\tilde \pi)$. Let $\AA'=\AA \setminus
\{B\}$. Let $A$ be the first in the top (last in the bottom) and
$E$ be the first in the bottom (last in the top) of $\tilde \pi$.
If $B=A$ let $\pi$ be obtained from $\tilde\pi$ by applying a top
arrow.
If $B=E$ let $\pi$ be obtained from $\tilde\pi$ by applying a
bottom arrow.
Otherwise, let $\pi=\tilde \pi$. The first two cases of this
definition are condensed in the following formula:
{\tiny
$$
\tilde\pi=\left(\begin{array}{cccc}
          A & \cdot & \cdot & E \\
          E & \cdot & \cdot & A \end{array}\right)
\mapsto
\pi=\left(\begin{array}{cccc}
          A & \cdot & \cdot & E \\
          E & A & \cdot & \cdot \end{array}\right)
\quad\text{or}\quad
\pi=\left(\begin{array}{cccc}
          A & E & \cdot & \cdot \\
          E & \cdot & \cdot & A \end{array}\right)
$$}
In all cases, the fact that $e_B\notin H(\tilde\pi)$ easily
implies that $e_B\notin H(\pi)$. Let $\pi'$ be obtained from $\pi$
by forgetting $B$. Then $\pi$ is a simple extension of $\pi'$ and,
from the equivalence of (1) and (3) in Lemma~\ref
{simplereduction}, we obtain that $g(\pi)=g(\pi')$.
\end{pf}

Lemma \ref{pinchingrauzy} says that replacing any permutation by
another in the same Rauzy class does not affect the strong
pinching property. By Lemma \ref{dg}, if a class is non-minimal
then it contains some permutation $\pi$ which is a simple
extension of some $\pi'$ with $g(\pi)=g(\pi')$. Then, Lemma \ref
{g=g'} says that the action of $\Pi(\pi')$ is conjugate to the
action of a submonoid of $\Pi(\pi)$. It follows that strong
pinching for $\Pi(\pi')$ implies strong pinching for $\Pi(\pi)$.
Repeating this procedure one eventually reaches a permutation in
some minimal class. This means that these results reduce the proof
of Theorem~\ref{strongpinching} to the case of minimal Rauzy
classes.

We will also need a special formulation of this reduction which
has a stronger conclusion, since it provides properties of a
specific submonoid.

\begin{lemma} \label {sex2}

Let $\pi' \in \ssigma(\AA')$ be such that $E \in \AA'$ is the last
in the top and the first in the bottom. Let $\AA''=\AA' \setminus
\{E\}$ and $\pi''$ be obtained from $\pi'$ by forgetting $E$.
Assume that $\pi''$ is irreducible and that $g(\pi')=g(\pi'')$.
Assume also that the action of $\Pi(\pi'')$ on $H(\pi'')$ twists
isotropic subspaces and is strongly pinching. Let $\tilde
\Pi(\pi') \subset \Pi(\pi')$ be the submonoid of all $\gamma$ such
that $\Theta(\gamma) \cdot e_E=e_E$. Then the action of $\tilde
\Pi(\pi')$ on $H(\pi')$ twists isotropic subspaces and is strongly
pinching.

\end{lemma}

\begin{pf}
Notice $\pi'$ is not a simple extension of $\pi''$, since $E$ is
the last in the top (recall Remark \ref{converse}), and so we can
not apply Lemma \ref{g=g'} directly. Let $d=\# \AA$. Let $D$ be
last in the bottom for $\pi'$. Then $D$ is in the $k$-th position
in the top row for some $1 \leq k \leq d-1$. In fact, $k\le d-2$
because we assume that $\pi''$ is irreducible. Let $\gamma'_0$ be
the bottom arrow starting at $\pi'$, and let $\tilde \pi'$ be the
end of $\gamma'_0$. Let $\gamma'_1$ be a sequence of $d-k-1$
bottoms starting at $\tilde \pi'$. Notice that $\gamma'_1$ ends at
$\pi'$. This is illustrated in the following formula:
{\tiny
$$
\pi'=\left(\begin{array}{cccccc} \cdot & \cdot & D & \cdot & \cdot & E \\
                                E & \cdot & \cdot & \cdot & \cdot & D \end{array}\right)
\quad{\overset {\gamma_0'} \to }\quad
\tilde\pi'=\left(\begin{array}{cccccc} \cdot & \cdot & D & E & \cdot & \cdot \\
                                E & \cdot & \cdot & \cdot & \cdot & D \end{array}\right)
\quad{\overset {\gamma_1'} \to } \pi' \quad
$$}
Then $\tilde \pi'$ is a simple extension of $\pi''$. By Lemma~\ref
{g=g'}, the action of $\LE_*(\Pi(\pi'')) \subset \Pi(\tilde \pi')$
on $H(\tilde \pi')$ is conjugate to the action of $\Pi(\pi'')$ on
$H(\pi'')$, and so it also twists isotropic subspaces and is
strongly pinching. Since $D \neq E$ is the winner of all arrows of
$\gamma'_0$ and $\gamma'_1$, $\Theta(\gamma_0')$ and
$\Theta(\gamma_1')$ preserve $e_E$. It follows that $\gamma'_0
\LE_*(\Pi(\pi'')) \gamma'_1$ is contained in the submonoid $\tilde
\Pi(\pi')$ that stabilizes $e_E$. We conclude as in the proof of
Lemma \ref {large} that $\tilde \Pi(\pi')$ twists isotropic
subspaces and is strongly pinching.
\end{pf}

%
%

\subsection{Proof of Theorem \ref {strongpinching}}

The proof is by induction on $\# \AA$. The case $\# \AA=2$ is easy
because any arrow, top or bottom, gives a parabolic element for the
action. Let us show that if strong pinching holds for $\# \AA=d-1 \geq 2$
then it holds for $\# \AA=d$.

Let $\pi \in \ssigma(\AA)$. As explained before, if $d>2g(\pi)$
then the result follows from the induction hypothesis using Lemmas
\ref {g=g'}, \ref {pinchingrauzy}, and \ref {dg}. So, we may
assume that $d=2g(\pi)$. Let $g=g(\pi)$. Notice that $g \geq 2$
because $d > 2$. By  Lemma \ref {pinchingrauzy}, the conclusion
does not change if we replace $\pi$ by any other permutation in
the same Rauzy class. By Lemma \ref {goodexists}, the Rauzy class
contains some good permutation $\tilde \pi \in \RR(\pi)$. Hence,
we may suppose from the start that $\pi$ is the permutation
obtained by applying a top arrow to $\tilde \pi$. Lemma \ref
{criterion} reduces the proof that the action of $\Pi(\pi)$ on
$H(\pi)=\R^\AA$ is strongly pinching to proving

\begin{lemma} \label{weareleft}

For every $C>0$ there exists $B \in \Theta(\Pi(\pi))$ for which $1$ is an
eigenvalue of geometric multiplicity 1 (the eigenspace has dimension $1$)
and
 \be \label{BAB}
 \theta_{g-1}(B)>0
 \quad\text{and}\quad
 \theta_k(B)>C\theta_{k+1}(B)
 \text{ for all\ \ } 1 \leq k \leq g-2.
 \ee

\end{lemma}

\begin{pf}

Let $A$ be first in the top (last in the bottom) and $E$ be first
in the bottom (last in the top) of $\pi$. Let $\pi'$ be obtained
by forgetting $A$ from $\pi$ and $\AA'=\AA \setminus \{A\}$. Then
$E$ is last in the top and the first in the bottom of $\pi'$, and
this ensures that $\pi'$ is irreducible. Moreover, $\pi$ is a
simple extension of $\pi'$, and $\pi'$ is a simple reduction of
both $\pi$ and $\tilde\pi$. This is illustrated in the formula
that follows:
{\tiny
$$
\tilde\pi=\left(\begin{array}{ccccc} A & \cdot & \cdot & E \\
                                     E & \cdot & \cdot & A\end{array}\right)
\quad\quad
\pi=\left(\begin{array}{ccccc} A & \cdot & \cdot & E \\
                               E & A & \cdot & \cdot \end{array}\right)
\quad\quad
\pi'=\left(\begin{array}{ccccc} \cdot & \cdot & E \\
                                E & \cdot & \cdot \end{array}\right)
$$}
Let $\gamma'_0 \in \Pi(\pi')$ be the sequence of $2g-2$ top arrows
starting (and ending) at $\pi'$, and let
$\gamma_0=\LE_*(\gamma'_0)$. Then $\gamma_0 \in \Pi(\pi)$ is a
sequence of $2g-1$ top arrows starting and ending at $\pi$. Notice
that $E$ is the winner in all these arrows, and so \be \label{U0}
\Theta(\gamma_0) =
\left(\begin{array}{ccccc} 1 & 0 & \cdots & 0 & 1 \\
                           0 & 1 & \cdots & 0 & 1 \\
                     \cdots & \cdots & \cdots & \cdots & \cdots \\
                           0 & 0 & \cdots & 1 & 1 \\
                           0 & 0 & \cdots & 0 & 1
\end{array}\right)
\ee where the $1$ column is the one indexed by $E$. More
precisely,
 \be \label {U01}
 \Theta(\gamma_0) \cdot e_x = e_x \quad\text{for all } x \in \AA
 \setminus \{E\},
 \ee
 \be \label {U02} \Theta(\gamma_0) \cdot e_E=\sum_{x \in \AA} e_x.
\ee Now let $\tilde \Pi(\pi') \subset \Pi(\pi')$ be the set of all
$\gamma'$ such that $\Theta(\gamma') \cdot e_E=e_E$.   We claim
that the action of $\tilde \Pi(\pi')$ on $H(\pi')$ twists
isotropic subspaces and is strongly pinching. This is easily
obtained from Lemma \ref{sex2}, as follows. Let $\pi''$ be
obtained from $\pi'$ by forgetting the letter $E$. Equivalently,
$\pi''$ is obtained from the good permutation $\tilde\pi$ by
forgetting $A$ and $E$. According to Lemma \ref{211},
$g(\tilde\pi)=g(\pi'')+1$. We also have $2g(\pi'') \le 2g(\pi')
\le \#\AA-1 = 2 g(\tilde\pi)-1$. Consequently, $g(\pi'')=g(\pi')$.
This means we are in a position to apply Lemma \ref{sex2}: since
the action of $\Pi(\pi'')$ on $H(\pi'')$ is strongly pinching and
twists isotropic subspaces, by the induction hypothesis and
Corollary \ref{twisting}, the same is true for the action of
$\tilde \Pi(\pi')$ on $H(\pi')$, as claimed. Then we may apply
Lemma \ref{more} to find $\gamma' \in \tilde \Pi(\pi')$ such that
(take $x=\gamma_0'$ and $x_0=\gamma'$ in the lemma) $1$ is not an
eigenvalue of $\gamma'$ acting on $H(\pi')$, and the Lyapunov
exponents of $\gamma'\gamma'_0$ acting on $H(\pi')$ satisfy \be
\label{(dois)} \theta_{g-1}(\gamma'\gamma'_0)>0
\quad\text{and}\quad \theta_i(\gamma'\gamma'_0)>C
\theta_{i+1}(\gamma'\gamma'_0) \quad\text{for all $1 \leq i \leq
g-2$.} \ee Let us show that $B=\Theta(\LE_*(\gamma'\gamma'_0))$
has the required properties.

Write $\gamma=\LE_*(\gamma')$, so that $B=\Theta(\gamma_0)\Theta(\gamma)$.
Since $E$ is never a winner for any of the arrows that form
$\gamma'$ (because $\gamma'$ is chosen in the stabilizer of
$e_E$), it is never a winner for $\gamma$ either. Moreover, noting
that $E$ is first and $A$ is second in the bottom of $\pi$, we see
that $A$ is never a winner nor a loser for $\gamma$. Together with
\eqref{U01} and \eqref{U02}, this gives
 \be \label {U1} \Theta(\gamma)
\cdot e_A = e_A \text{ and } \Theta(\gamma) \cdot e_E=e_E
 \ee
 \be
\label{U2} \langle \Theta(\gamma) \cdot e_x,e_A \rangle=0
\quad\text{for all } x \in \AA \setminus \{A\}.
 \ee
In other words, the matrix of $\Theta(\gamma)$ has the form \be
\label{U}
\Theta(\gamma)=\left(\begin{array}{ccccc}1 & 0 & \cdot & 0 & 0 \\
                          0 & * & \cdot & * & 0 \\
                          \cdot & \cdot & \cdot & \cdot & \cdot \\
                          0 & * & \cdot & * & 0 \\
                          0 & * & \cdot & * & 1 \end{array}\right)
\ee This implies that $B=\Theta(\gamma_0)\Theta(\gamma)$ satisfies
\be \label {BeA} B \cdot e_A=e_A, \ee \be \label {BeE} \langle B
\cdot e_x,e_E \rangle=\langle B \cdot e_x,e_A \rangle
\quad\text{for all } x \in \AA \setminus \{A\}. \ee By \eqref
{BeA}, we have that $1$ is an eigenvalue of $B$. Let us check that
its geometric multiplicity is $1$. Indeed, otherwise there would
exist $h \in \R^\AA \setminus \{0\}$ such that $B \cdot h=h$ and
$\langle h,e_A \rangle=0$.  By \eqref {BeE}, \be \langle h, e_E
\rangle = \langle B \cdot h,e_E \rangle = \langle B \cdot h,e_A
\rangle = \langle h,e_A \rangle=0. \ee Let $P:\R^\AA\to\R^{\AA'}$
be the natural projection and $h'=P \cdot h$. Notice that $P^*
\cdot h'=h$, because $h$ is orthogonal to $e_A$. In particular,
$h'$ is non-zero. We have
 \be \label {PP}
 \Theta(\gamma_0') = P \Theta(\gamma_0) P^*, \quad \Theta(\gamma')
 = P \Theta(\gamma) P^*, \quad \Theta(\gamma_0')\Theta(\gamma') = P
 B P^*.
 \ee
This implies \be \label {hker1} \Theta(\gamma'_0)\Theta(\gamma')
\cdot h' = P B P^* \cdot h' = P B \cdot h = P \cdot h=h'. \ee We
also have $\langle h',e_E \rangle = \langle h, e_E \rangle = 0$.
Now, \eqref {U01} and \eqref{PP} imply that $\Theta(\gamma'_0)
\cdot e_x = e_x$ for all $x \in \AA' \setminus \{E\}$.
Consequently, \be \label{hker2} \Theta(\gamma'_0) \cdot h' = h'
\quad\text{or, equivalently,}\quad h'=\Theta(\gamma'_0)^{-1} \cdot
h'. \ee From \eqref{hker1} and \eqref{hker2} we immediately get
that \be \label{hker3} \Theta(\gamma') \cdot h' = h'. \ee But
\eqref {U1} implies that $\Theta(\gamma') \cdot e_E =
P\Theta(\gamma) P^* \cdot e_E=e_E$, and so \eqref{hker3} implies
that $1$ is an eigenvalue of $\Theta(\gamma')$ with geometric
multiplicity $2$ at least.  Since $H(\pi')$ has codimension $1$ in
$\R^{\AA'}$, this implies that $1$ is an eigenvalue of
$\Theta(\gamma')$ acting on $H(\pi')$, which contradicts the
choice of $\gamma'$. This proves that $1$ is an eigenvalue of $B$
with geometric multiplicity $1$, as claimed. To obtain the
properties in \eqref {BAB}, notice that \eqref {BeA} and \eqref
{PP} imply that the matrix of $B$ has a block triangular form,
with the matrix of $\Theta(\gamma'_0)\Theta(\gamma')$ as a
diagonal block. It follows that the eigenvalues of $B$ are,
precisely, the eigenvalues of $\Theta(\gamma'_0)\Theta(\gamma')$,
together with an additional eigenvalue $1$ associated to the
eigenvector $e_A$. Observe that this eigenvalue must have {\it
algebraic} multiplicity $2$, since the action is symplectic.
Therefore, the Lyapunov spectrum of $B$ consists of the Lyapunov
spectrum of $\Theta(\gamma'_0)\Theta(\gamma')$ on $H(\pi')$
together with two additional zero Lyapunov exponents. Thus,
\eqref{(dois)} implies \eqref {BAB}. This finishes the proof of
Lemma \ref{weareleft}.
\end{pf}

At this point the proof of Theorem \ref {strongpinching} is complete.
\qed

\begin{rem}

In the above argument we concatenate two Rauzy paths to generate another one
exhibiting some parabolic behavior. It would be interesting to find a
geometric interpretation, in terms of the Teichm\"uller flow, of how
that behavior arises.

\end{rem}

\begin{rem}

Some of the richness properties we have been discussing, like
twisting, depend only on the group generated by the monoid
(because invariance under the monoid is equivalent to invariance
under the group generated), while others, such as strong pinching,
depend on the monoid itself. It is interesting to investigate how
large is the generated group. The only obvious restriction is that
this group preserves the integer lattice $H(\pi) \cap \Z^\AA$.
What is its Zariski closure (which plays a role in \cite {GM})~?
Such questions had already been raised by Zorich in \cite {Z4},
Appendix A.3, where he made some specific conjectures. We believe
it is possible to proceed further in this direction with the
arguments of the present paper, particularly our induction
procedure in terms of relations between Rauzy classes.

\end{rem}

\section{Simplicity of the spectrum} \label {simpl}

\subsection{A criterion for simplicity of the spectrum} \label {statement of
criterion}

Let $(\Delta,\mu_0)$ be a probability space, and let $T:\Delta \to
\Delta$ be a transformation such that there exists a finite or
countable partition of $\Delta$ into sets $\Delta^{(l)}$, $l \in
\Lambda$ of positive measure such that $T:\Delta^{(l)} \to \Delta$
is an invertible measurable transformation and
$T_*(\mu_0|\Delta^{(l)})$ is equivalent to $\mu_0$.  Let $\Omega$
be the set of finite sequences of elements of $\Lambda$, including
the empty sequence. If $\l=(l_1,\ldots ,l_m) \in \Omega$, let
$$
\Delta^\l=\{x \in \Delta: \, T^k(x) \in \Delta^{(l_{k+1})} \text{ for } 0 \leq k<m\},
$$
and $T^\l=T^m:\Delta^\l \to \Delta$. The $\Delta^\l$ have positive
measure and $T^\l$ is an invertible measurable transformation. We
say that $(T,\mu_0)$ has {\it approximate product structure} if
there exists $C>0$ such that
 \be \label{stexp}
 \frac {1} {C}
 \leq \frac {1} {\mu_0(\Delta^\l)} \frac {dT_*^\l(\mu_0|\Delta^\l)}{d\mu_0}
 \leq C \quad\text{for all } \l \in \Omega.
 \ee
Notice that $T$ is measurable with respect to the {\it cylinder
$\sigma$-algebra} of $T$, that is, the $\sigma$-algebra generated
by the $\Delta^\l$, $\l \in \Omega$. It is not difficult to check,
using \eqref{stexp}, that $T$ is ergodic with respect to the
cylinder $\sigma$-algebra, and there is a unique probability
measure $\mu$ on the cylinder $\sigma$-algebra which is invariant
under $T$ and is absolutely continuous with respect to $\mu_0$.
Moreover,
 \be
 \frac {1} {C} \leq \frac {d\mu} {d\mu_0} \leq C
 \ee
and $(T,\mu)$ has approximate product structure as well.

Let $(T,\mu)$ have approximate product structure and $H$ be some
finite dimensional vector space. Let $A^{(l)} \in \SL(H)$, $l \in
\Lambda$, and define $A:\Delta \to \SL(H)$ by $A(x)=A^{(l)}$ if $x
\in \Delta^{(l)}$. We will say that $(T,A)$ is a {\it locally
constant cocycle.} The {\it supporting monoid} of $(T,A)$ is the
monoid generated by the $A^{(l)}$, $l \in \Lambda$. By ergodicity
and the Oseledets theorem~\cite{O}, if the cocycle is measurable,
that is, if $\int \ln \|A(x)\| d\mu(x)<\infty$, then it has a well
defined Lyapunov spectrum.

\begin{main} \label {simplicity criterion}
Let $(T,A)$ be a locally constant measurable cocycle.  If the supporting
monoid is simple then the Lyapunov spectrum is simple.
\end{main}

This is an adaptation of the main result in \cite{AV} to the
present situation. For completeness, a proof is provided in
Appendix~\ref{simplicity criterion}.

\subsection{Locally constant projective cocycles}

Fix any $p \geq 2$.  We call $\P \R^p_+$ the {\it standard
simplex}. A {\it projective contraction} is a projective
transformation taking the standard simplex into itself or, in
other words, it is the projectivization of some matrix $B \in
\GL(p,\R)$ with non-negative entries. The image of the standard
simplex by a projective contraction is called a {\it simplex.}

A {\it projective expanding map} $T$ is a map $\cup \Delta^{(l)}
\to \Delta$, where $\Delta$ is a simplex compactly contained in
the standard simplex, the $\Delta^{(l)}$ form a finite or
countable family of pairwise disjoint simplexes contained in
$\Delta$ and covering almost all of $\Delta$, and
$T^{(l)}=T:\Delta^{(l)} \to \Delta$ is a bijection such that
$(T^{(l)})^{-1}$ is the restriction of a projective contraction.

\begin{lemma} \label {projective1}

If $T:\cup \Delta^{(l)} \to \Delta$ is a projective expanding map
then it has approximate product structure with respect to Lebesgue
measure.

\end{lemma}

\begin{pf}

This follows from the well-known fact that the logarithm of the
Jacobian of a projective contraction relative to Lebesgue measure
is (uniformly) Lipschitz with respect to the projective metric.
See the proof of Lemma 2.1 of \cite {AF} or the Appendix of
\cite{AV} for details.
\end{pf}

%
%
%
%

\subsection{Proof of Theorem \ref{main}}

Let $\gamma \in \Pi(\RR)$ be any path such that $\P \Delta_\gamma$
is compactly contained in $\P \Delta^1_\RR$. This is easy to
satisfy: for every $x \in \P \Delta_\RR$, the connected component
of $x$ in $\P \Delta^n_\RR$ is compactly contained in $\P
\Delta^1_\RR$ for every $n$ sufficiently large. Let $T$ be the
first return map to $\P \Delta_\gamma$ under the Zorich map. We
define a linear skew-product $(x,h) \mapsto (T(x),A(x) \cdot h)$
over $T:\P \Delta_\gamma \to \P \Delta_\gamma$, by setting
$A(x)=B^Z(R_Z^{m-1}(x)) \cdots B^Z(x)$, where $m=m(x)$ is the
first return time of $x$ to $\P \Delta_\gamma$. The map $T$
preserves the probability measure
$$
\mu=\frac {1} {\mu_Z(\P \Delta_\gamma)} \mu_Z|(\P \Delta_\gamma),
$$
and the Lyapunov exponents of $(T,A)$ are obtained multiplying the
Lyapunov exponents of the cocycle $(R_Z,B^Z)$ by ${1}/{\mu_Z(\P
\Delta_\gamma)}$. Notice that $T$ is a projective expanding map.
So, by Lemma~\ref {projective1}, $(T,A)$ is a locally constant
cocycle.

Let $\pi'$ be the start of $\gamma$. Let $\gamma_0$ be an arrow
ending at $\pi'$ and which is either top or bottom according to
whether $\gamma$ starts by a bottom or a top arrow. Let $\pi$ be
the start of $\gamma_0$. Let $\gamma_1$ be a path starting by
$\gamma$ and ending at $\pi$. Then the monoid $\BB$ generated by
the $A(x)$ contains $\Theta(\gamma_1 \Pi(\pi) \gamma_0)$. Hence
$\BB$ is a submonoid of $\Theta(\Pi(\pi'))$ containing
$\Theta(\gamma_1 \gamma_0 \Pi(\pi') \gamma_1 \gamma_0)$.  Since,
the action of $\Theta(\Pi(\pi'))$ on $H(\pi')$ is simple, Lemma
\ref{large} gives that the monoid $\BB$ generated by the $A(x)$ is
also simple. Now Theorem~\ref{main} is a consequence of
Theorem~\ref {simplicity criterion}. \qed

\appendix

\section{Proof of Theorem~\ref{simplicity criterion}}

Here we prove our sufficient criterion for simplicity of the Lyapunov
spectrum. See also \cite{AV}.

\subsection{Symbolic dynamics}

Let us say that two locally constant cocycles $(T,A)$ and
$(T',A')$ are equivalent if $\mu(\Delta^\l)=\mu'(\Delta'^\l)$ for
all $\l \in \Omega$ and $A^{(l)}=A'^{(l)}$ for all $l \in
\Lambda$. It is clear that two measurable locally constant
cocycles have the same Lyapunov spectrum if they are equivalent.
Any locally constant cocycle is equivalent to a {\it symbolic}
locally constant cocycle, that is, one for which the base
transformation is the shift $(x_i)_{i \in \N} \mapsto (x_{i+1})_{i
\in \N}$ over $\Lambda^{\N}$, where $\N=\{n \in \Z,\, n \geq 0\}$.
Thus, it is enough to prove Theorem~\ref {simplicity criterion} in
the symbolic case.  Let us introduce some convenient notation to
deal with this case.

Let $\Sigma=\Lambda^{\N}$. Given $\l=(l_1,\ldots ,l_m) \in
\Omega$, let $\Sigma^\l$ be the set of $(x_i)_{i=0}^\infty$ such
that $x_i=l_{i+1}$, $0 \leq i<m$. We also write $|\l|=m$. Let
$f:\Sigma \to \Sigma$ be the shift map. A cylinder on $\Sigma$ is
a set of the form $f^{-n}(\Sigma^\l)$ with $n \in \N$.  Denote by
$f^\l$ the restriction of $f^m$ to $\Sigma^\l$. A probability
measure on $\Sigma$ is defined once specified on cylinders. Thus,
an $f$-invariant probability measure is completely determined by
its value on the $\Sigma^\l$. Given an $f$-invariant probability
$\mu$, let $\phi_\mu:\Omega \to \R^+$ be defined by
$\phi_\mu(\l)=\mu(\Sigma^\l)$. For $\l_1=(l_1,\ldots ,l_m)$ and
$\l_2=(l_{m+1},\ldots ,l_{m+n})$ we denote $\l_1 \l_2=(l_1,\ldots
,l_{m+n}$).

We will say that $\mu$ has {\it bounded distortion} for $f$ if it
gives positive measure to all cylinders and there exists
$C(\mu)>0$ such that
 \be
 \frac {1} {C(\mu)}
 \leq \frac {\phi_\mu(\l_1 \l_2)} {\phi_\mu(\l_1)\phi_\mu(\l_2)}
 \leq C(\mu), \quad\text{for all } \l_1,\l_2 \in \Omega.
 \ee
Notice that this last condition is equivalent to
 \be
 \frac {1} {C(\mu)}
 \leq \frac {1} {\mu(\Sigma^\l)} \frac {df_*^\l(\mu|\Sigma^\l)}{d\mu}
 \leq C(\mu), \quad\text{for all } \l \in \Omega.
 \ee
Hence, $\mu$ has bounded distortion for $f$ if and only if
$(f,\mu)$ has approximate product structure.

\subsection{The inverse limit}

Let $\Sigma_-=\Lambda^{\Z \setminus \N}$ and $\widehat
\Sigma=\Lambda^\Z$.  Then, $\widehat \Sigma=\Sigma_- \times
\Sigma$. Let $\pi_-$ and $\pi$ be the coordinate projections. Let
$\hat f:\widehat \Sigma \to \widehat \Sigma$ be the shift map.
Then $\pi \circ \hat f=f \circ \pi$. Let $\widehat
\Sigma^\l=\Sigma_- \times \Sigma^\l$ and $\widehat
\Sigma^\l_-=\hat f^{m}(\Sigma^\l)$, with $m=|\l|$, and
$\Sigma^\l_- \subset \Sigma_-$ be such that $\widehat
\Sigma^\l_-=\Sigma^\l_- \times \Sigma$.  A cylinder on $\widehat
\Sigma$ is a set of the form $\hat f^n(\widehat \Sigma^\l)$, $n
\in \Z$.

A probability measure on $\widehat \Sigma$ is defined once
specified on cylinders.  An $\hat f$-invariant probability measure
is, thus, completely determined by its values on the $\widehat
\Sigma^\l$. Thus, there is a bijection $\hat \mu \mapsto \mu$
between $\hat f$- and $f$-invariant probability measures given by
$\pi_* \hat \mu=\mu$.  We call $\hat \mu$ the {\it lift} of $\mu$.
We say that $\hat \mu$ has bounded distortion for $\hat f$ if
$\mu$ has bounded distortion for $f$. This implies that $\hat f^n$
is ergodic with respect to $\hat \mu$ for every $n \geq 1$. We
also denote by $\mu_-$ the projection of $\hat\mu$ to $\Sigma_-$.

%
%
%
%
%

\subsection{Invariant section}

Fix some probability measure $\mu$ with bounded distortion and let
$\hat \mu$ be its lift. Given a locally constant measurable
cocycle $(f,A)$, we may consider its lift $(\hat f,\widehat A)$,
defined by $\widehat A(x) = A(\pi(x))$. It naturally acts on
$\widehat \Sigma \times \grass(k,H)$, for any $1 \leq k \leq \dim
H-1$. Let $k$ be fixed. Given $\l \in \Omega$, let $\xi_\l=A^\l
\cdot E^+_k(A^\l)$ if defined, otherwise choose $\xi_\l$
arbitrarily. For $x \in \Sigma_-$ and $n\in\N$, let $\l(x,n) \in
\Omega$ be such that $|\l(x,n)|=n$ and $x \in \Sigma^{\l(x,n)}_-$.

\begin{thm} \label {xi}

Let $(f,A)$ be a symbolic locally constant measurable cocycle, and
assume the supporting monoid is simple. Then,
\begin{enumerate}
\item There exists a measurable function $\xi:\Sigma_- \to
\grass(k,H)$ such that $\hat \xi=\xi \circ \pi_-$ satisfies $\hat
\xi(\hat f(x))=\hat A(x) \cdot \hat \xi(x)$,

\item For $\mu_-$-almost all $x \in \Sigma_-$ we have $\ln
\sigma_k(A^{\l(x,n)})-\ln \sigma_{k+1}(A^{\l(x,n)}) \to \infty$
and $\xi_{\l(x,n)} \to \xi(x)$.

\item For every hyperplane section $S \subset \grass(k,H)$ there
exists a positive $\mu_-$-measure set of $x \in \Sigma_-$ such
that $\xi(x) \notin S$.
\end{enumerate}

\end{thm}

The proof of this theorem will take several steps.

\subsection{Measures on Grassmannians}

Let $\hat m$ be a probability on $\widehat\Sigma \times
\grass(k,H)$. We say that $\hat m$ is a {\it $u$-state} if its
projection on $\widehat \Sigma$ is $\hat \mu$ and there exists
$C(\hat m)>0$ such that for any Borel set $X\subset \grass(k,H)$
and any $\l_0,\l,\l' \in \Omega$, we have \be \frac {1}
{\mu(\Sigma^{\l'})} \hat m(\Sigma^{\l_0}_- \times \Sigma^{\l'}
\times X) \leq C(\hat m) \frac {1} {\mu(\Sigma^\l)} \hat
m(\Sigma^{\l_0}_- \times \Sigma^\l \times X). \ee
It is easy to give examples of  $u$-states: take any probability
measure $\nu$ in $\grass(k,H)$ and let $\hat m=\hat \mu \times
\nu$ (in this case $C(\hat m)=C(\mu)^2$ works). From this we can
get examples of $u$-states invariant under $(\hat f,\hat A)$, as
follows.


Let $\NN$ be the space of probability measures on $\widehat
\Sigma\times\grass(k,H)$ that project to $\hat \mu$ on
$\widehat\Sigma$. We introduce on $\NN$ the smallest topology such
that the map $\eta \mapsto \int \psi d\eta$ is continuous, for
every bounded continuous function $\psi:\widehat \Sigma \times
\grass(k,H) \to \R$. We will call this the weak$^*$ topology.
Notice that $\NN$ is a compact separable space. This is easy to
see in the following alternative description of the topology. Let
$K_n \subset \widehat \Sigma$, $n\ge 1$ be disjoint compact sets
such that $\hat\mu(K_n)>0$ and $\sum \hat \mu(K_n)=1$.  Let
$\NN_n$ be the space of measures on $K_n \times \grass(k,H)$ that
project to $(\hat \mu \mid K_n)$.  The usual weak$^*$ topology
makes $\NN_n$ a compact separable space. Given $\eta \in \NN$, let
$\eta_n \in \NN_n$ be obtained by restriction of $\eta$.  This
identifies $\NN$ with $\prod \NN_n$, and the weak$^*$ topology on
$\NN$ corresponds to the product topology on $\prod \NN_n$.


\begin{lemma}

Let $\hat m_0$ be a $u$-state.  For every $n \geq 0$,
\be
\hat m^{(n)}_0=(\hat f,\widehat A)^n_* (\hat m_0)
\ee
is a  $u$-state and $C(\hat m^{(n)}_0) \leq C(\hat m_0) C(\mu)^2$.

\end{lemma}

\begin{pf}

We must show for $\l_0,\l,\l' \in \Omega$ and $X \subset \grass(k,H)$
measurable that
 \be
\frac {\mu(\Sigma^\l)} {\mu(\Sigma^{\l'})} \hat m^{(n)}_0
(\Sigma^{\l_0}_- \times \Sigma^{\l'} \times X)
 \leq C(\hat m_0) C(\mu)^2
\hat m^{(n)}_0(\Sigma^{\l_0}_- \times \Sigma^\l \times X).
 \ee
It is enough to consider the case when $|\l_0| \geq n$.  In this case write
$\l_0=\l_1 \l_2$ with $|\l_2|=n$.  Then
 \begin{align}
 \frac{\mu(\Sigma^\l)}{\mu(\Sigma^{\l'})}
 \hat m^{(n)}_0(\Sigma^{\l_0}_- \times \Sigma^{\l'} \times X)
 & =
 \frac{\mu(\Sigma^\l)}{\mu(\Sigma^{\l'})}
 \hat m_0(\Sigma^{\l_1}_- \times \Sigma^{\l_2 \l'} \times (A^{\l_2})^{-1} \cdot X)
 \\
 \nonumber & \leq C(\hat m_0)
 \frac{\mu(\Sigma^\l)}{\mu(\Sigma^{\l'})}
 \frac {\mu(\Sigma^{\l_2 \l'})} {\mu(\Sigma^{\l_2 \l})}
 \hat m_0(\Sigma^{\l_1}_- \times \Sigma^{\l_2 \l} \times (A^{\l_2})^{-1} \cdot X)
 \\
 \nonumber & \leq
 C(\hat m_0) C(\mu)^2
  \hat m_0(\Sigma^{\l_1}_- \times \Sigma^{\l_2 \l} \times (A^{\l_2})^{-1} \cdot X)
 \\
 \nonumber & =
 C(\hat m_0) C(\mu)^2 \hat m^{(n)}_0(\Sigma^{\l_0}_- \times \Sigma^\l \times X),
\end{align}
as claimed.
\end{pf}

For any fixed $C>0$, the space of $u$-states $\hat m$ with $C(\hat
m)\le C$ is a convex compact subset in the weak$^*$ topology. So,
the next statement implies the existence of invariant $u$-states.

\begin{cor} \label {existence}

Let $\hat m_0$ be a $u$-state.  Let $\hat m$ be a Cesaro weak$^*$
limit of the sequence
\be
(\hat f,\widehat A)^n_* (\hat m_0).
\ee
Then $\hat m$ is an invariant $u$-state with $C(\hat m) \leq C(\hat m_0)
C(\mu)^2$.

\end{cor}

\begin{lemma} \label {conditionals}

Let $\hat m$ be a probability measure on $\widehat \Sigma \times
\grass(k,H)$.  For $x \in \Sigma_-$, let $\hat m^{(n)}(x)$ be the
probability on $\grass(k,H)$ obtained by normalized projection of
$\hat m$ restricted to
$\Sigma^{\l(x,n)}_-\times\Sigma\times\grass(k,H)$. Then, for
$\mu_--$almost every $x$, the sequence $\hat m^{(n)}(x)$ converges
in the weak$^*$ topology to some probability measure $\hat m(x)$.

\end{lemma}

\begin{pf}

Since $\grass(k,H)$ is a compact metric space, and keeping in mind
the definition of the weak$^*$ topology, it is enough to show that
if $\phi:\grass(k,H) \to \R$ is continuous then for almost every
$x \in \Sigma_-$ the integral $\phi^{(n)}(x)$ of $\phi$ with
respect to $\hat m^{(n)}(x)$ converges as $n\to\infty$. This is a
simple application of the martingale convergence theorem. A direct
proof goes as follows.

We may assume that $0 \leq \phi(z) \leq 1$, $z \in \grass(k,H)$.
Consider the measure $\nu \leq \mu_-$ defined by
$$
\nu(X)=\int_{X \times \Sigma \times \grass(k,H)} \phi(z) d\hat
m(x,y,z).
$$
Then, $\nu(\Sigma_-^{\l(x,n)})=\phi^{(n)}(x)
\mu_-(\Sigma_-^{\l(x,n)})$. Let $\phi^+(x)=\limsup \phi^{(n)}(x)$
and $\phi^-(x)=\liminf \phi^{(n)}(x)$.  We have to show that for
every $a<b$ the set $U_{a,b}=\{x \in \Sigma_-,\,
\phi^-(x)<a<b<\phi^+(x)\}$ has zero $\mu_--$measure.  To do this
it is enough to show that $\nu(U_{a,b}) \geq
b\mu_-(U_{a,b})-\epsilon$ and $\nu(U_{a,b}) \leq a
\mu_-(U_{a,b})+\epsilon$ for every $\epsilon>0$. We will show only
the first inequality, the second one being analogous. Fix an open
set $U\supset U_{a,b}$ with $\mu_-(U \setminus U_{a,b})<\epsilon$.
Let
$$
N(x)=\min \{n,\, \Sigma^{\l(x,n)}_- \subset U\} \quad\text{and}\quad
n(x)=\min \{n \geq N(x),\, \phi^{(n)}(x)>b\}.
$$
Let $U^N=\cup_{x \in U_{a,b}} \Sigma^{\l(x,n(x))}_-$. Then
$U_{a,b} \subset U^N \subset U$ and there exists a (finite or
countable) sequence $x_j \in U_{a,b}$ such that $U^N=\cup_j
\Sigma^{\l(x_j,n(x_j))}_-$. In particular, $\nu(U^N)-\nu
(U_{a,b})<\epsilon$ and
$$
\nu(U^N)=\sum_j \nu(\Sigma^{\l(x_j,n(x_j))}_-) \geq \sum_j b
\mu_-(\Sigma^{\l(x_j,n(x_j))}_-) \geq b \mu_-(U_{a,b}).
$$
This implies that $\nu(U_{a,b}) \geq b \mu_-(U_{a,b})-\epsilon$,
as required.
\end{pf}

%
%
%
%
%

\begin{lemma}\label{previous}

There exists $N>0$, $\delta>0$, $m>0$, $\l_p$, $\l_q$ and
$\l_{q_i}$, $1 \leq i \leq m$ with the following properties:

\begin{enumerate}

\item $|\l_p|=|\l_q|=|\l_{q_i}|=N$,

\item The matrix $A^{\l_p}$ has simple Lyapunov spectrum,

\item $A^{\l_q} \cdot F \cap G=\{0\}$ for every $F \in
\grass(k,H)$ and $G\in\grass(\dim H-k,H)$ which are sums of
eigenspaces of $A^{\l_p}$,

\item For every $F \in \grass(k,H)$ and $G \in \grass(\dim H-k,H)$
such that $F$ is a sum of eigenspaces of $A^{\l_p}$, there exists
$i$ such that the $A^{\l_{q_i}} \cdot F$ and $G$ form angle at
least $\delta$.

\end{enumerate}

\end{lemma}

\begin{pf}

Since $\BB$ is simple, there exists $\l_{\tilde p}$ such that the
matrix $A^{\l_{\tilde p}}$ has simple Lyapunov spectrum.  Let
$\BB' \subset \BB$ be the monoid consisting of the $A^\l$ where
$|\l|$ is a multiple of $|\l_{\tilde p}|$. Then $\BB'$ is a large
submonoid of $\BB$, in the sense of Lemma~\ref {large}, and so
$\BB'$ is simple. By the definition of twisting and Lemma~\ref
{twist}, there exists $\l_{\tilde q}$ such that $|\l_{\tilde q}|$
is a multiple of $|\l_{\tilde p}|$ and $A^{\l_{\tilde q}} \cdot F
\cap G=\{0\}$ for every $F \in \grass(k,H)$ and every $G\in
\grass(\dim H-k,H)$ which are sums of eigenspaces of
$A^{\l_{\tilde p}}$. For the same reasons, for every
$G\in\grass(\dim H-k,H)$ there exists $\l_{\tilde q(G)}$ such that
$|\l_{\tilde q(G)}|$ is a multiple of $|\l_{\tilde p}|$ and
$A^{\l_{\tilde q(G)}} \cdot F \cap G=\{0\}$ for every $F \in
\grass(k,H)$ which is a sum of eigenspaces of $A^{\l_{\tilde p}}$.
By compactness of $\grass(\dim H-k,H)$, one can  choose finitely
many $\l_{\tilde q_i}$ among all the $\l_{\tilde q(G)}$, so that
for every $G \in \grass(k,H)$,
$$
\max_i\min_F\operatorname{angle}\big(A^{\l_{\tilde q_i}}\cdot F, G) > 0,
$$
where the minimum is over all $F \in \grass(k,H)$ that are sums of
eigenspaces of $A^{\l_{\tilde p}}$. Using compactness once more,
the expression on the left is even uniformly bounded below by some
$\delta>0$. Let $N$ be the maximum of $|\l_{\tilde p}|$,
$|\l_{\tilde q}|$ and $\l_{\tilde q_i}$.  Let $\l_p$, $\l_q$ and
$\l_{q_i}$ be obtained from $\l_{\tilde p}$, $\l_{\tilde q}$ and
$\l_{\tilde q_i}$ by adding at the beginning as many copies of
$\l_{\tilde p}$ as necessary to get to length $N$.  One easily
checks the required properties.
\end{pf}


For $\l \in \Omega$, let $\l^n$ be obtained by repeating $\l$ exactly
$n$ times. Given a probability measure $\rho$ on $\grass(k,H)$, let
$\rho^\l$ be the push-forward of $\rho$ under
$A^\l:\grass(k,H) \to \grass(k,H)$.

%
%
%

\begin{lemma}

Given $\epsilon>0$ and any probability measure $\rho$ in
$\grass(k,H)$, there is $n_0=n_0(\epsilon,\rho)$ and, given any
$\l_0 \in \Omega$, there is $i=i(\l_0)$ such that
$\rho^\l(B)>1-\epsilon$ for every $n>n_0$, where $\l=\l^n_p \l_q
\l^n_p \l_{q_i} \l_0$ and $B$ is the $\epsilon$-neighborhood of
$\xi_\l$.

\end{lemma}

\begin{pf}

By (2) in Lemma \ref{previous}, when $n$ is large most of the mass
of $\rho^{\l^n_p}$ is concentrated near sums of eigenspaces of
$A^{\l_p}$ of dimension $k$. Then most of the mass of
$\rho^{\l^n_p \l_q}$ is concentrated near the $A^{\l_q}$-image of
those sums. Using (3) in Lemma \ref{previous}, it follows that,
when $n$ is large most of the mass of $\rho^{\l^n_p \l_q \l^n_p}$
is concentrated near the subspace $F\in\grass(k,H)$ given by the
sum of the eigenspaces associated to the $k$ eigenvalues of
$A^{\l_p}$ with largest norms. Now let $G(\l_0)\in\grass(\dim H
-k, H)$ be spanned by eigenvectors of $(A^{\l_0})^*A^{\l_0}$
corresponding to $\dim H-k$ smallest singular values; if
$\sigma_k(A^{\l_0})>\sigma_{k+1}(A^{\l_0})$ then the only choice
is $G(\l_0)=E^-_k(A^{\l_0})$. Then, in general, the family
$A^{\l_0}$, $\l_0$ is equicontinuous restricted to the subset of
$k$-dimensional subspaces whose angle to $G(\l_0)$ is larger than
any fixed $\delta/2$. Choose $\delta>0$ as in Lemma
\ref{previous}, and let $i=i(\l_0)$ be as in (4) of that lemma,
for this choice of $F$ and $G=G(\l_0)$. Then most of the mass of
$\rho^{\l^n_p \l_q \l^n_p \l_{q_i}}$ is concentrated near
$A^{\l_{q_i}}\cdot F$, and so most of the mass of $\rho^\l$ is
concentrated near $A^{\l_{q_i} \l_0} \cdot F$, as long as $n$ is
large; moreover, the equicontinuity property allows us to take the
condition on $n$ uniform on $\l_0$. In particular, we get that
$\rho^{\l}$ converges to the Dirac measure at $A^{\l_{q_i} \l_0}
\cdot F$ when $n$ goes to infinity.

Now it suffices to show that $\xi_{\l}$ is close to $A^{\l_{q_i}
\l_0} \cdot F$ when $n$ is large, uniformly on $\l_0$. This can be
done applying the previous argument to the case when $\rho$ is the
Lebesgue measure on $\grass(k,H)$, and using (2) in Lemma~\ref
{non-trivial} for $x_n=A^{\l^n_p \l_q \l^n_p}$: we conclude that
$\ln \sigma_k(A^\l)-\ln \sigma_{k+1}(A^\l) \to \infty$ and that
$\xi_{\l}=A^{\l} \cdot E^+_k(A^{\l})$ converges to $A^{\l_{q_i}
\l_0} \cdot F$ when $n\to\infty$, as claimed.
\end{pf}

\begin{lemma} \label {x_n}

There exists $C_0=C_0(f,A,\mu)>0$ and a sequence of sets $X_n
\subset \Omega$ such that $|\l|\ge n$ for all $\l\in X_n$, the
$\widehat \Sigma^\l_-$, $\l \in X_n$ are pairwise disjoint with
$\sum_{\l \in X_n} \mu_-(\Sigma^\l_-)>C^{-1}_0$ and, given any
probability measure $\rho$ on $\grass(k,H)$ and any
$n>n_0(\epsilon,\rho)$, we have
 \be \rho^\l(B)>(1-\epsilon)
\quad\text{for all}\quad \l \in X_n
 \ee
where $B$ is an $\epsilon$-neighborhood of $\xi_\l$.

\end{lemma}

\begin{pf}

Let us order the elements of $\l \in \Omega$ by inclusion of the
$\widehat \Sigma^\l_-$.
Let $X^0_n$ be the collection of elements of $\Omega$ of the form
$\l^n_p \l_q \l^n_p \l_{q_i} \l_0$ where $|\l_0|$ is a multiple of
$(2n+2) N$.  Let $X^1_n \subset X^0_n$ be the maximal
elements.  Then the $\widehat \Sigma^\l_-$, $\l \in X^1_n$ are disjoint and
(by ergodicity of $\hat f^{(2n+2) N}$)
\be
\sum_{\l \in X^1_n} \mu_-(\Sigma^\l_-)=1.
\ee
Notice that if $l^n_p \l_q \l^n_p \l_{q_i} \l_0 \in X^1_n$ for some $1 \leq
i \leq m$ then $l^n_p \l_q \l^n_p \l_{q_i} \l_0 \in X^1_n$ for every
$1 \leq i \leq m$.
Let $X_n$ be the set of $\l \in X^1_n$ such that $\l=l^n_p \l_q \l^n_p
\l_{q_i} \l_0$ with $i=i(\l_0)$ as in the previous lemma.  Then
\be
\sum_{\l \in X_n} \mu_-(\Sigma^\l_-) \geq C(\mu)^{-2} \min_{1 \leq i \leq m}
\mu_-(\Sigma^{\l_{q_i}}_-).
\ee
To conclude, apply the previous lemma.
\end{pf}

\noindent{\it Proof of Theorem~\ref {xi}.} Let $\hat m$ be an
invariant $u$-state, given by Corollary~\ref {existence}, and let
$\nu$ be its projection on $\grass(k,H)$. For almost every $x \in
\Sigma_-$, let $\hat m^{(n)}(x)$ and $\hat m(x)$ be as in
Lemma~\ref {conditionals}. Notice that we have, for any Borel set
$Y \subset \grass(k,H)$, \be \frac {1} {C(\hat m)}
 \leq \frac {\hat m^{(n)}(x)(Y)} {\nu^{\l(x,n)}(Y)} \leq C(\hat m).
\ee In particular, $\hat m(x)$ is equivalent to any limit of
$\nu^{\l(x,n)}$.  By Lemma~\ref {x_n}, there is $Z \subset
\Sigma_-$ with $\mu_-(Z) \geq C_0^{-1}$ such that $\nu^{\l(x,n)}$
accumulates at a Dirac mass for all $x\in Z$.  Using
Lemma~\ref{conditionals}, it follows that $\hat m(x)$ is a Dirac
mass for almost every $x \in Z$. By ergodicity and equivariance,
$\hat m(x)$ is a Dirac mass for almost every $x \in \Sigma_-$. We
will denote the support of this Dirac mass  by $\xi(x)$. Hence,
$\xi(x)$ is the support of $\lim \nu^{\l(x,n)}$. Let $\hat \xi=\xi
\circ \pi_-$. Then $\hat \xi(f(x,y))=\hat A(x,y)\cdot\hat
\xi(x,y)$. Notice that the push-forward of $\mu_-$ by $\xi$ is
equal to $\nu$. In particular, the support of $\nu$ is forward
invariant under the action of the supporting monoid $\BB$. Observe
that the support of $\nu$ can not be contained in any hyperplane
section $S$: otherwise, $\cap_{x \in \BB} x^{-1} \cdot S$ would be
a non-trivial invariant linear arrangement, and this can not exist
because $\BB$ is simple. It follows from part (2) of Lemma~\ref
{non-trivial} that $\ln\sigma_k(A^{\l(x,n)})-\ln
\sigma_{k+1}(A^{\l(x,n)}) \to \infty$ and $\xi(x)=\lim_{n \to
\infty} \xi_{\l(x,n)}$. \qed

\begin{rem}

It follows that the push-forward of $\hat \mu$ by $(x,y,z) \mapsto
(x,y,\xi(x))$ is an invariant $u$-state. It is not difficult to
see that it must coincide with $\hat m$. Thus, the invariant
$u$-state is unique. See also Remark~5.4 in \cite{BV}.

\end{rem}

\subsection{The inverse cocycle}

Given $\l \in \Omega$, let $\xi^-_\l=E^-_k(A^\l)$ if defined,
otherwise choose $\xi^-_\l$ arbitrarily. For $x \in \Sigma$, let
$\l^-(x,n) \in \Omega$ be such that $|\l^-(x,n)|=n$ and $x \in
\Sigma^{\l^-(x,n)}$.

\begin{thm} \label {xi^-}

Let $(f,A)$ be a symbolic locally constant measurable cocycle, with simple
supporting monoid.
\begin{enumerate}
\item There exists a
measurable function $\xi^-:\Sigma \to \grass(\dim H-k,H)$ such that
$\hat \xi^-=\xi^- \circ \pi$ satisfies $\hat \xi^-(\hat f(x))=
\hat A(x) \cdot \hat \xi^-(x)$,

\item For $\mu-$almost every $x \in \Sigma$, $\ln \sigma_k
(A^{\l^-(x,n)})-\ln \sigma_{k+1}(A^{\l^-(x,n)}) \to \infty$ and
$\xi^-_{\l^-(x,n)} \to \xi^-(x)$.

\item For every hyperplane section $S \subset \grass(\dim H -
k,H)$ there exists a positive $\mu-$measure set of $x \in \Sigma$
such that $\xi^-(x) \notin S$.
\end{enumerate}

\end{thm}

\begin{pf}

The cocycle $(\hat f,\widehat A)^{-1}$ is also measurable. Let
$J:\widehat \Sigma \to \widehat \Sigma$ be given by $J (x_i)_{i
\in \Z}=(x_{-i-1})_{i \in \Z}$. Then $J$ conjugates $(\hat
f,\widehat A)^{-1}$ to a locally constant symbolic cocycle $(\hat
f,\hat B)$, where $\hat B$ is defined by $\hat B(x)=\widehat
A^{-1}(J(x))=\widehat A(\hat f^{-1}(J(x)))^{-1}$. The
corresponding supporting monoid is also simple, by Lemma~\ref
{b-1}. The result then follows by application of Theorem~\ref {xi}
to $(\hat f,\hat B)$, with $\dim H-k$ in the place of $k$. The
invariant sections are related by $\hat\xi^-_{\hat
A}=\hat\xi_{\hat B} \circ J$.
\end{pf}

\subsection{Proof of Theorem~\ref {simplicity criterion}}

Let $1 \leq k \leq \dim H-1$.  We must show that
$\theta_k(f,A)>\theta_{k+1}(f,A)$. Let $\xi:\Sigma_- \to
\grass(k,H)$ be as in Theorem~\ref {xi}, and let $\xi^-:\Sigma \to
\grass(\dim H-k,H)$ be as in Theorem~\ref {xi^-}. We claim that
$\xi(x) \cap \xi^-(y)=\{0\}$ for $\hat\mu-$almost every $(x,y) \in
\Sigma_- \times \Sigma$. Indeed, if this was not the case then, by
ergodicity, we would have $\xi(x) \cap \xi^-(y) \neq \{0\}$ for
$\hat\mu-$almost every $(x,y)$. By Fubini theorem and bounded
distortion, it would follow that $\xi(x) \cap \xi^-(y) \neq \{0\}$
for $\mu_--$almost every $x \in \Sigma^-$ and $\mu-$almost every
$y \in \Sigma$. This would imply that $\xi(x)$ is contained in the
hyperplane section dual to $\xi^-(y)$ for $\mu_--$almost every $x
\in \Sigma^-$, contradicting Theorem~\ref{xi}.

Let $F^u \in \grass(k,H)$ and $F^s \in \grass(\dim H-k,H)$ be
subspaces in the support of $\xi_* \mu_-$ and $\xi^-_* \mu$,
respectively, such that $F^u \cap F^s=\{0\}$. Choose an inner
product so that $F^u$ and $F^s$ are orthogonal. Let $\epsilon_0>0$
be such that if $F$ and $F'$ belong to the balls of radius
$\epsilon_0$ around $F^u$ and $F^s$, respectively, then $F \cap
F'=\{0\}$. Fix $0<\epsilon_3<\epsilon_2<\epsilon_1<\epsilon_0$.
Take $0<\epsilon_4<\epsilon_3$ small and let $\widehat X=X_-
\times X$ be the set of $(x,y)$ such that $\xi(x)$ belongs to the
ball of radius $\epsilon_4$ around $F^u$ and $\xi^-(y)$ belongs to
the ball of radius $\epsilon_4$ around $F^s$. Let $n(x,y)\ge 0$ be
minimum such that $A^{\l(x,n)}$ takes the ball of radius
$\epsilon_1$ around $F^u$ into the ball of radius $\epsilon_3$
around $F^u$ whenever $n>n(x,y)$ is such that $f^{-n}(x,y) \in
\widehat X$. Then $n(x,y)<\infty$ for almost every $(x,y) \in X_-
\times X$, by Theorems \ref {xi} and \ref {xi^-} and Lemma~\ref
{compactlimit}.

Let $\widehat Z \subset \widehat X$ be a positive measure set such
that the minimum of the first return time $r(x,y)$ of $f$ to
$\widehat Z$ is bigger than the maximum of $n(x,y)$ over $(x,y)
\in \widehat Z$. Let $R$ denote the return map to $\widehat Z$.
For every $(x,y) \in \widehat Z$, choose $C(x,y)\in\SL(H)$ in a
measurable fashion in the $c^{-1} \epsilon_4$-neighborhood of
$\id$, such that $C(x,y) \cdot \xi(x)=F^u$ and $C(x,y) \cdot
\xi^-(y)=F^s$. Let $B(x,y)=C(R(x,y)) A^\l C(x,y)^{-1}$, where
$\l=\l^-(y,r(x,y))$. Then the Lyapunov exponents of $(R,B)$
relative to the normalized restriction of $\hat\mu$ to $\widehat
Z$ are obtained by multiplying the Lyapunov exponents of $(f,A)$
by ${1}/{\hat \mu(\widehat Z)}$. Taking $\epsilon_4$ sufficiently
small, we guarantee that $B(x,y)$ takes the closure of the ball of
radius $\epsilon_2$ around $F^u$ into its interior. It follows
from Lemma~\ref {fthetaargument} that
$\sigma_k(B(x,y))>\sigma_{k+1}(B(x,y))$ for almost every $(x,y)
\in \widehat Z$, and $E^+_k(B(x,y))=F^u$ and $E^-_k(B(x,y))=F^s$.
We conclude that \be \theta_k(R,B)-\theta_{k+1}(R,B)
  \geq \frac {1} {\hat \mu(Z)}
  \int_Z \Big[\ln \sigma_k(B(x,y))-\ln \sigma_{k+1}(B(x,y))\Big]\,
   d \hat \mu(x,y)>0,
\ee
which implies the conclusion of the theorem.
\qed

\subsection{Generalizations}

Theorem~\ref{simplicity criterion} may be seen as a criterion for
simplicity of the Lyapunov spectrum for products of matrices
generated by a stationary stochastic process. Although we have
stated it in the case when there is a finite or countable number
of states, our method can actually be used to prove the following
more general result. Let $(T:(\Delta,\mu) \to (\Delta,\mu),
A:\Delta \to \SL(H))$ be a measurable cocycle. Given $k\ge 1$ and
any set $Z \in \SL(H)^k$, we denote
 \be
Z^\wedge = \{x\in\Delta : (A(x), \ldots, A(T^{k-1}(x))) \in Z\}.
 \ee
We say that the cocycle has \emph{approximate local product
structure} if there exists $C \geq 1$ such that, for any $m, n\ge
1$ and measurable sets $X \in \SL(H)^n$, $Y \in \SL(H)^m$,
 \be
 C^{-1} \mu (X^\wedge) \mu(Y^\wedge)
 \le \mu\ ((X\times Y)^\wedge)
 \le C \mu (X^\wedge) \mu(Y^\wedge).
 \ee

\begin{thm}

Let $(T:(\Delta,\mu) \to (\Delta,\mu), A:\Delta \to \SL(H))$ be a
measurable cocycle with approximate local product structure, and
assume that the submonoid of $\SL(H)$ generated by the support of
the measure $A_*\mu$ is simple. Then the Lyapunov spectrum of
$(T,A)$ is simple.

\end{thm}


\begin{thebibliography}{MMY}

\bibitem[AF]{AF} Avila, A.; Forni, G.
Weak mixing for interval exchange transformations and translation
flows. Preprint (www.arxiv.org) 2004, to appear Annals of Math.

\bibitem[AV1]{AV0}
Avila, A.; Viana, M. Dynamics in the moduli space of Abelian
differentials. Preprint (www.preprint.impa.br) 2005.

\bibitem[AV2]{AV}
Avila, A.; Viana, M. Simplicity of Lyapunov spectra: a sufficient
criterion. In preparation.

\bibitem[BGV]{BGMV}
Bonatti, C.; Gomez-Mont, X.; Viana, M.
G\'en\'ericit\'e d'exposants de Lyapunov non-nuls pour des
produits d\'eterministes de matrices.
Ann. Inst. H. Poincare Anal. Non Lineaire 20 (2003), no. 4, 579--624.

\bibitem[BV]{BV}
Bonatti, C.; Viana, M.
Lyapunov exponents with multiplicity $1$ for deterministic products of
matrices.
Ergod. Th. \& Dynam. Syst. 24 (2004), 1295--1330.

\bibitem[Fo]{F}
Forni, G. Deviation of ergodic averages for area-preserving flows
on surfaces of higher genus. Ann. of Math. 155 (2002), no. 1,
1--103.

\bibitem[Fu]{Fu}
Furstenberg, H.
Noncommuting random products.
Trans. Amer. Math. Soc. 108 1963 377--428.

\bibitem[GM]{GM}
Gol'dsheid, I. Ya.; Margulis, G. A.
Lyapunov exponents of a product of random matrices.
Uspekhi Mat. Nauk 44 (1989), no. 5(269), 13--60; translation in
Russian Math. Surveys 44 (1989), no. 5, 11--71.

\bibitem[GR]{GR}
Guivarc'h, Y.; Raugi, A.
Products of random matrices: convergence theorems.
Random matrices and their applications (Brunswick, Maine, 1984), 31--54,
Contemp. Math., 50,
Amer. Math. Soc., Providence, RI, 1986.

\comm{
Guivarc'h, Y.; Raugi, A.
   Proprietes de contraction d'un semi-groupe de matrices inversibles.
   Coefficients de Liapunoff d'un produit de matrices aleatoires
   independantes.
Israel J. Math. 65 (1989), no. 2, 165--196.
}

\bibitem[Ko]{K}
Kontsevich, M.
Lyapunov exponents and Hodge theory.
The mathematical beauty of physics (Saclay, 1996), 318--332,
Adv. Ser. Math. Phys., 24, World Sci. Publishing, River Edge, NJ, 1997.

\bibitem[KZ]{KZ}
Kontsevich, M.; Zorich, A.
Connected components of the moduli spaces of Abelian differentials with
prescribed singularities. Invent. Math. 153 (2003), no. 3, 631--678.

\bibitem[Le]{L}
Ledrappier, F.
Positivity of the exponent for stationary sequences of matrices.
Lyapunov exponents (Bremen, 1984), 56--73,
Lecture Notes in Math., 1186,
Springer, Berlin, 1986.

\bibitem[MMY]{MMY}
Marmi, S.; Moussa, P.; Yoccoz, J.-C. The cohomological equation
for Roth type interval exchange transformations. Preprint Scuola
Normale Superiore di Pisa, 2003.

\bibitem[Ma]{M}
Masur, H. Interval exchange transformations and measured
foliations. Ann. of Math. 115 (1982), no. 1, 169--200.

\bibitem[Os]{O}
Oseledets, V. I. A multiplicative ergodic theorem: {L}yapunov
characteristic numbers for dynamical systems. Trans. Moscow Math.
Soc. 19 (1968), 197--231.

\bibitem[Ra]{R}
Rauzy, G.
Echanges d'intervalles et transformations induites.
Acta Arith. 34, (1979), no. 4, 315--328.

\bibitem[Sc]{S}
Schwartzman, S. Asymptotic cycles. Ann. of Math. 66 (1957),
270--284.

\bibitem[Ve1]{V1}
Veech, W. Gauss measures for transformations on the space of
interval exchange maps. Ann. of Math. 115 (1982), no. 1, 201--242.

\bibitem[Ve2]{V2}
Veech, W. The Teichm\"uller geodesic flow. Ann. of Math. 124
(1986), no. 3, 441--530.

\bibitem[Vi]{V}
Viana, M. Almost all cocycles over any hyperbolic system have
non-vanishing Lyapunov exponents. Preprint (www.preprint.impa.br)
2005.

\bibitem[Zo1]{Z1}
Zorich, A.
Asymptotic flag of an orientable measured foliation on a surface.
Geometric study of foliations (Tokyo, 1993), 479--498,
World Sci. Publishing, River Edge, NJ, 1994.

\bibitem[Zo2]{Z2}
Zorich, A. Finite Gauss measure on the space of interval exchange
transformations. Lyapunov exponents. Ann. Inst. Fourier 46 (1996),
no. 2, 325--370.

\bibitem[Zo3]{Z3}
Zorich, A.
Deviation for interval exchange transformations.
Ergodic Theory Dynam. Systems 17 (1997), no. 6, 1477--1499.

\bibitem[Zo4]{Z4}
Zorich, A.
How do the leaves of a closed $1$-form wind around a surface?
Pseudoperiodic topology, 135--178,
Amer. Math. Soc. Transl. Ser. 2, 197, Amer. Math. Soc.,
Providence, RI, 1999.

\end{thebibliography}
\end{document}